\documentclass[twocolumn]{autart}                                                           
\usepackage{latexsym,amssymb,amsmath, amsbsy,amsopn, amstext,multicol,multirow}
\usepackage{graphicx,epsfig,tikz,pgf,hyperref,cancel,color,float,ifpdf}
\usepackage{extarrows}
\usepackage{amsmath, amssymb,stfloats}
\usepackage{setspace}
\usepackage{tikz, pgfplots,subcaption}
 
\graphicspath{{./}}
\newcommand{\ud}{\mathrm{d}}
\def\e{\epsilon}
\def\l{\lambda}
\def\d{\delta}

\def\b{\beta}
\def\m{\mu}

\def\S{\Sigma}
\def\b{\beta}

\def\R{\mathbb{R}}
\def\Z{\mathbb{Z}}

\def\ra{\rightarrow}
\def\lra{\longrightarrow}

\allowdisplaybreaks[4]

\newtheorem{definition}{\bfseries Definition}
\newtheorem{example}{\bfseries Example}
\newtheorem{assumption}{\bfseries Assumption}
\newtheorem{theorem}{\bfseries Theorem}

\newtheorem{corollary}{\bfseries Corollary}
\newtheorem{lemma}{\bfseries Lemma}

\newtheorem{remark}{\bfseries Remark}

\begin{document}
\begin{frontmatter}
\title{Passivity-Based Analysis of Sampled and Quantized Control Implementations}
\thanks[footnoteinfo]{The material in this paper was partially presented at the 54th IEEE Conference on Decision and Control, Dec 15-18, 2015,  Osaka, Japan. Corresponding author: Xiangru Xu.}
\author[XXU]{Xiangru Xu}\ead{xiangru.xu@wisc.edu},   
\author[NOZ]{Necmiye Ozay}\ead{necmiye@umich.edu},
\author[VGU]{Vijay Gupta}\ead{vgupta2@nd.edu}

\address[XXU]{Department of Mechanical Engineering, University of Wisconsin, Madison, WI, 53706, USA} 
\address[NOZ]{Department of Electrical Engineering \& Computer Science, University of Michigan, Ann Arbor, MI, 48109, USA} 
\address[VGU]{Department of Electrical Engineering, University of Notre Dame, Notre Dame, IN, 46556, USA} 


\begin{abstract}
This paper studies the performance of a continuous controller when implemented on digital devices via sampling and quantization, by leveraging passivity analysis. 
Degradation of passivity indices from a continuous-time control system to its sampled, input and output quantized model is studied using a notion of quasi-passivity. Based on that,  the passivity property of a feedback-connected system where the continuous controller is replaced by its sampled and quantized model is studied, and conditions that ensure the state boundedness of the interconnected system are provided. Additionally, the approximate bisimulation-based control implementation where the controller is replaced by its approximate bisimilar symbolic model whose states are also quantized is analyzed. Several examples are provided to illustrate the theoretical results.	

\end{abstract}
\end{frontmatter}

\section{Introduction}\label{sec:introduction}
Cyber-physical control systems consist of software-based controllers interacting with physical processes (i.e., the plant to be controlled). While the control design can be done using continuous-time and continuous-space methods, in order to guarantee the desirable operation of the closed-loop system, it is important to keep in mind the restrictions the software implementation imposes. In particular, sampling and quantization are prevalent in real control implementations. Therefore, sampled-data and quantized control have attracted the attention of researchers for decades. However, many problems are still open, especially when dynamics of the system are nonlinear or when both sampling and quantization are considered \cite{chen2012optimal,brockett2000quantized,elia2001stabilization,fu2005sector}. This work presents a unified framework for analyzing control implementations using tools from passivity theory. In addition to analyzing the effects of sampling and quantization of the input and output signals within this framework, we also consider quantization of the internal states of the controller. This enables analysis of controllers implemented with finitely many bits, and can be relevant in extremely resource-constrained settings like micro-robotics or applications where a computer that can do floating-point arithmetic is not feasible, which constitute the motivation for our work.

The concept of passivity emerged from the study of energy dissipation in circuit analysis, which, roughly speaking, means that a system can not generate ``internal energy'' on its own. In the paper \cite{JCWillems72}, Willems systematized the theory of dissipativity (with passivity as its special case) using the concepts of storage function and supply rate, where dissipativity is defined as the property that the rate of increase of the storage function is not larger than the supply rate. Since then, the fundamental connections between dissipativity, stability (either in the sense of Lyapunov or input/output) and optimality for control systems have been explored, making passivity/dissipativity theory a widely-used tool in control theory \cite{BrogDissBook07,Khalil2000,sepulchre2012constructive}. One important property of passive systems is the \emph{compositionality} - the parallel or feedback connection of two passive systems remains passive. This  property offers an effective method for analyzing the stability of large-scale, interconnected systems in a compositional way 
\cite{lee2014passivity,wen2004unifying}. 
Several challenging problems in cyber-physical systems, such as quantization, time delay
and packet drops, which are induced by the digital devices or communication networks, have also been studied under the passivity/dissipativity-based framework  (see ~\cite{hirche2012human,stramigioli2005sampled,gao2007passivity,wang14,antsaklis2013control} and references therein). 


A general form of dissipativity that has been  well-studied is the QSR-dissipativity whose supply rate is a quadratic function of the input and output. A special case of the QSR-dissipativity allows one to use  two real numbers, which are termed \emph{passivity indices} (or passivity levels), to characterize the shortage/excess of passivity for a system. 
These two numbers can be either positive or negative, where a negative value implies a shortage of passivity while a positive value implies an excess. Passivity indices are rather useful in the analysis of nonlinear and interconnected systems; for instance,  \cite{hill1976stability,hill1977stability} showed that the feedback interconnection of two non-passive systems is finite gain $\mathcal{L}_2$ stable provided that the shortage  of passivity of one component can be compensated for by the excess  of passivity of another component,  \cite{madeira2016asymptotic} provided a controller design method for asymptotic stabilization of a class of nonlinear systems using passivity indices, and \cite{ZhuFeedbackACC14} and \cite{zhu2017passivity} presented conditions on passivity indices of a closed-loop system in the continuous-time setting and with the event-triggered mechanism, respectively.

It is known that passivity is not preserved under sampling or quantization in general (e.g., see \cite{SenPos02,antsaklis2013control,OishiCDC10}). 
However, passivity indices, as a quantitative ``abstraction'' of the passivity property of a system, can be used to study how the properties of two related systems (e.g., a continuous-time system and its sampled or quantized system) differ from each other. 
Some interesting results have been obtained along this direction;
for example, \cite{OishiCDC10} analyzed the passivity degradation for a system under time sampling with the zero-order hold and the ideal sampler, \cite{xia2017approx}  showed that passivity indices of one system can be inferred from an approximate model of that system when the norm of the model error is small and the approximate model has an excess of passivity. Despite these interesting results, passivity analysis for systems that are both sampled and quantized, as well as its control implementations in a feedback loop, are still lacking, to the best of our knowledge.


Symbolic models for dynamical systems provide a unified framework for studying the interactions of software and physical phenomena \cite{PauloBook09}.
Such symbolic models are sampled and quantized (in state, input and output), which are used in control software synthesis from high-level specifications \cite{PolaAuto08}. 
The basic workflow in these approaches is to first compute an approximate symbolic model of the plant based on (bi)simulation relations, then synthesize a discrete controller for the symbolic plant model, and finally refine that controller and compose it with the plant.  The constructed controller can be implemented in software and the closed-loop system is guaranteed to satisfy the desired high-level specifications. A major limitation of these techniques is the curse of dimensionality: the complexity of the symbolic plant model 
grows exponentially with the dimension of the state-space of the plant (when quantized uniformly).
To tackle this computational problem, compositional approaches have been utilized. For instance,  \cite{majumdar2016compositional} proposed a compositional symbolic abstraction method for a class of continuous-time nonlinear control systems using the notion of approximate disturbance bisimulation and discussed the related controller synthesis problem, \cite{dallal2015compositional} provided a compositional control design method, inspired by small gain theorem and assume-guarantee reasoning, for the feedback composition of two discrete systems with a persistency specification. Passivity, which possesses the nice compositionality property, has also been used to tame the computational complexity of symbolic control. For instance, by using simulation functions and the dissipativity-based approach,  \cite{rungger2015compositional,zamani2017compositional} investigated compositional abstractions of a network of continuous-time control systems to a lower dimensional, continuous-time interconnected system. However, results that combine the bisimulation-based finite state abstraction and the passivity/dissipativity-based approach are still rare.



In this paper, we analyze the passivity property of sampled and quantized controller implementations with passivity indices as the main tool. Specifically, the question we consider is - suppose that a continuous, dynamic controller has been designed to ensure specified passivity indices for the closed-loop system, what guarantees on passivity can be inherited if a sampled and quantized controller  is implemented, and under what conditions the closed-loop system is (ultimately) bounded? 
To this end, we first propose the notion of quasi-passivity and the strong detectability for discrete-time control systems, and give several lemmas that relate ultimate boundedness, strong detectability and passivity indices.
Then, we explore how passivity degrades from a continuous-time system to its sampled, input and output quantized model, where a set of results quantifying the passivity degradation are derived, which relate the passivity indices, the sampling time and the quantization precision. Based on that, we study the implementation of the continuous controller by replacing it with its sampled, input and output quantized model, and we provide conditions (on the passivity indices and the strong detectability of the interconnected components) under which the state of the closed-loop system eventually converges into a compact set whose size can be made arbitrarily small by choosing the quantization precision small enough. Finally, we consider implementing an approximate bisimulation-based symbolic controller whose states are also quantized, in the closed-loop system, and give conditions that guarantee the ultimate boundedness of the closed-loop system.

A preliminary version of the paper appeared in the conference publication  \cite{xupasscdc15}. In this paper, we extend the results in \cite{xupasscdc15} in the following important ways: {\color{black} Lemma 3, 4 and 5 that relate quasi-passivity, strong detectability and state boundedness are added, Theorem 1, 2 and 3 that present passivity degradation results are derived using a relaxed assumption (i.e., Assumption \ref{assgain}), Theorem 4 and 5 that present state boundedness results for the closed-loop system are new, all the complete proofs are included, and a few new examples are provided. }




\emph{Notation and definitions.}
The set of non-negative integers  and real numbers are denoted as $\Z_{\geq 0}$ and $\R_{\geq 0}$, respectively. 
The $\ell_2$ and $\ell_\infty$ norm of a vector $x\in \R^n$ are denoted as $|x|$ and $|x|_\infty$, respectively;  the $\ell_\infty$ norm of a function $\phi:\Z_{\geq 0}\rightarrow \R^n$ is denoted as $\|\phi\|:=\sup_{k\in\Z_{\geq 0}}|\phi[k]|$.  
For any $A\subseteq \R^n$ and $\mu>0$, $[A]_\m:=\{a\in A|a_i=k_i\m,k_i\in \Z,i=1,2,...,n\}$.
A relation $R\subset A\times B$ is identified with the map $R:A\ra 2^B$, which is defined by $b\in R(a)$ if and only if $(a,b)\in R$. For a set $S$, the set $R(S)$ is defined as $R(S)=\{b\in B: \exists a\in S, (a,b)\in R\}$. Given a relation $R\subset A\times B$, $R^{-1}$ denotes the inverse relation of $R$, i.e., $R^{-1}:=\{(b,a)\in B\times A: (a,b)\in R\}$. 
A continuous function $f: \R_{\geq 0}\rightarrow \R_{\geq 0}$ belongs to class $\mathcal{K}$ if it is strictly increasing and $f(0)=0$, and $f$ belongs to class $\mathcal{K}_{\infty}$ if $f\in\mathcal{K}$ and $f(r)\rightarrow\infty$ as $r\rightarrow \infty$; a continuous function $f: \R_{\geq 0}\times \R_{\geq 0}\rightarrow \R_{\geq 0}$ belongs to class $\mathcal{KL}$ if for each fixed $s$, function $f(r,s)\in \mathcal{K}_{\infty}$ with respect to $r$ and for each each fixed $r$, function $f(r,s)$ is decreasing with respect to $s$ and $f(r,s)\ra 0$ as $s\ra \infty$. 




\section{Preliminaries and Problem Statement}\label{sec:preliminary}

\subsection{Transition Systems and Approximate Bisimulation}\label{sec:contsystem}

A continuous-time control system is a tuple $\S=(X,U,Y,f,h)$ where $X\subseteq \R^n$ is a set of states, $U\subseteq\R^m$ is a set of inputs, $Y\subseteq \R^s$ is a set of outputs, $f:X\times U\ra \R^n$ is Lipschitz continuous, and $h:X \times U\ra \R^s$ is continuous. The state, input and output of $\S$ at time $t\in\R_{\geq 0}$ are denoted by $x(t)$, $u(t)$, $y(t)$, respectively, and their evolution is governed by:
\begin{equation}\label{ctd}
\begin{aligned}
\dot{x}(t) &=f(x(t),u(t)),\\
y(t) & =h(x(t),u(t)), \;\; \forall t \in \R_{\geq 0}.
\end{aligned}
\end{equation}
We assume that 
$f(0,0)=0$ and $h(0,0)=0$. 
We also assume that given any sufficiently regular control input signal ${\bf u}:[0,T]\rightarrow U$ with $T> 0$ and any initial condition $x_0\in X$, there exist a unique state trajectory $\bf{x}$ and a corresponding output trajectory $\bf{y}$ defined on $[0,T]$ satisfying $x(0)=x_0$ and Eq.~\eqref{ctd}. Denote by $\textbf{x}(\tau,x_0,\textbf{u})$  the state reached at time $\tau$ under the input $\textbf{u}$ from the initial state $x_0$ of $\S$.


A discrete-time control system is a tuple $\S_d=(X,U,Y,f_d,h_d)$ where
$X\subseteq \R^n$ is a set of states, $U\subseteq\R^m$ is a set of inputs, $Y\subseteq \R^s$ is a set of outputs, $f_d:X\times U\ra X$ and $h_d:X \times U\ra \R^s$ are both continuous maps.
The state, input and output of $\S_d$ at time step $k\in \Z_{\geq 0}$ are denoted by $x[k]$, $u[k]$, $y[k]$, respectively, and their evolution is governed by:
\begin{equation}\label{dtd}
\begin{aligned}
x[k+1] &=f_d(x[k],u[k]),\\
y[k] & =h_d(x[k],u[k]), \;\; \forall k \in \Z_{\geq 0}.
\end{aligned}
\end{equation}
The state and output trajectories of the system $\S_d$ are discrete-time signals satisfying Eq.~\eqref{dtd}. 




\begin{definition}\label{FTS}
	A transition system is a quintuple $T=(Q,L,O,\lra,H)$, where:
	\begin{itemize}
		\item $Q$ is a set of states;
		\item $L$ is a set of inputs;
		\item $O$ is a set of outputs;
		\item $\lra\subset Q\times L\times Q$ is the transition relation;
		\item $H:Q\times L\lra O$ is the output function.
	\end{itemize}
\end{definition}

In the following, we denote an element $(q,\ell,p)\in\lra$ in a transition relation by $q\xlongrightarrow{\ell}p$ where $p,q\in Q$, $\ell\in L$. 

Given a continuous-time system $\Sigma=(\R^n,U,\R^s,f,h)$ and a sampling time $\tau$, we suppose that the control inputs of $\S$ are piecewise constant, that is, $u(t)=u((k-1)\tau)$ for any $t\in [(k-1)\tau,k\tau),k\in \Z_{\geq 0}$. Then, similar to \cite{PolaAuto08}, we define a transition system 
%
$T_\tau(\Sigma)=(X_1,U_1,Y_1,\xlongrightarrow[1]{},H_1)$ associated with the time-sampling of $\S$ as follows:
\begin{itemize}
	\item $X_1=\R^n$;
	\item $U_1=U$;
	\item $Y_1=\R^s$;
	\item $p\xlongrightarrow[1]{u} q$ if $\textbf{x}(\tau,p,\textbf{u})=q$ where $\textbf{u}:[0,\tau)\rightarrow \{u\}$, $u\in U_1$;
	\item $H_1(x,u)=h(x,u)$.
\end{itemize}

We interpret the trajectories of $T_\tau(\Sigma)$ in discrete-time, that is, it has an equivalent representation in terms of a discrete-time control system as in \eqref{dtd}, where its state, input and output at time step $k\in \Z_{\geq 0}$ are denoted by $x[k]$, $u[k]$, $y[k]$, respectively. Note that $T_\tau(\Sigma)$ can be obtained by putting $\Sigma$ between  a zero-order hold device (D/A) and an uniform sampler (A/D).

By further quantizing the state and input spaces of $T_\tau(\Sigma)$, we obtain an {\color{black}infinitely countable} transition system\footnote{Note that the sets of state, input and output of $T_{\tau\mu\eta}(\Sigma)$ are all countable.} $T_{\tau\mu\eta}(\Sigma)=(X_2,U_2,Y_2,\xlongrightarrow[2]{},H_2)$ for some $\tau,\mu,\eta>0$ as follows  (see also \cite{PolaAuto08}):
\begin{itemize}
	\item $X_2=[\R^n]_\eta$;
	\item $U_2=[U]_\mu$;
	\item $Y_2=\R^s$;
	\item $p\xlongrightarrow[2]{u} q$ if $|\textbf{x}(\tau,p,\textbf{u})-q|_\infty\leq \eta/2$ where $\textbf{u}:[0,\tau)\rightarrow \{u\}$, $u\in U_2$;
	\item $H_2(x,u)=h(x,u)$.
\end{itemize}


Bisimulation is a binary relation between two transition systems, which, roughly speaking, requires the two systems match each other's behavior \cite{PauloBook09}. In \cite{GirardTAC07},
the exact bisimulation was generalized to  $\e$-approximate bisimulation, which allows the states of two transition systems to be within certain bounds. To further capture the input and output behaviors of transition systems, we  consider the following \emph{$(\e,\mu)$-approximate bisimulation} relation adopted from \cite{julius2009approximate}. 


\begin{definition}\label{dfnemubi}
	Given two transition systems $T_1=(Q_1,L,O,\xlongrightarrow[1]{},H_1)$ and $T_2=(Q_2,L,O,\xlongrightarrow[2]{},H_2)$  {\color{black}where $Q_1\subseteq Q_2$, the sets $Q_1,Q_2$ are equipped with the same metric $\textbf{d}_s$ defined for $Q_2$}, and the input set $L$ is equipped with the metric $\textbf{d}_l$, for any $\e,\mu\in\R^+$, a relation $R\subset Q_1\times Q_2$ is said to be an $(\e,\mu)$-approximate bisimulation relation between $T_1$ and $T_2$, if for any $(q_1,q_2)\in R$:
	\begin{enumerate}
		\item[(i)] $\textbf{d}_s(q_1,q_2)\leq \e$;
		\item[(ii)] $q_1\xlongrightarrow[1]{\ell_1}p_1$ implies the existence of $\ell_2\in L$ such that $\textbf{d}_l(\ell_1,\ell_2)\leq \mu$, $q_2\xlongrightarrow[2]{\ell_2}p_2$ and $(p_1,p_2)\in R$;
		\item[(iii)] $q_2\xlongrightarrow[2]{\ell_2}p_2$ implies the existence of
		$\ell_1\in L$ such that $\textbf{d}_l(\ell_1,\ell_2)\leq \mu$, $q_1\xlongrightarrow[1]{\ell_1}p_1$ and $(p_1,p_2)\in R$.
	\end{enumerate}
\end{definition}

If there exists an $(\e,\mu)$-approximate bisimulation relation $R$ between $T_1$ and $T_2$ such that $R(Q_1)=Q_2$ and $R^{-1}(Q_2)=Q_1$, $T_1$ is said to be $(\e,\mu)$-bisimilar to $T_2$, which is denoted as $T_1\cong^{(\e,\mu)}T_2$.


\begin{remark}
	{\color{black} The $(\e,\mu)$-approximate bisimulation in Def. \ref{dfnemubi} is different from that in \cite{julius2009approximate}: the finite transition system in \cite{julius2009approximate} has an observation map $\langle\cdot\rangle: Q\rightarrow O$ but without the output function $H$, where $Q_1\subseteq Q_2$ is not required either. In fact, we can consider the  observation map of $\S$, $T_{\tau}(\Sigma)$ or $T_{\tau\mu\eta}(\Sigma)$ above as an identity mapping, which is the case in \cite{PolaAuto08}. The output function  $H$ defined in this paper is used for the passivity analysis.}
\end{remark}

\begin{definition}\cite{angeli2002lyapunov}
	{\color{black}The continuous-time system $\S$ in \eqref{ctd} is called \emph{incrementally input-to-state stable} ($\d$-ISS)  if it is forward complete and there exist functions $\b_1\in \mathcal{KL}$ and $\b_2\in \mathcal{K}_{\infty}$ such that for any $t\in\R_{\geq 0}$, any initial state $x_1,x_2\in \R^n$ and any input $\textbf{u},\textbf{v}$, it holds that 
		$|\textbf{x}(t,x_1,\textbf{u})\!-\!\textbf{x}(t,x_2,\textbf{v})|\leq\!\b_1(|x_1\!-\!x_2|,t)\!+\!\b_2(\|\textbf{u}\!-\!\textbf{v}\|).$}
\end{definition}

\begin{lemma}\label{thmsi}
	Consider the continuous-time control system $\S$ in \eqref{ctd} and any desired precision $\e>0$. If $\S$ is $\d$-ISS satisfying $|\textbf{x}(t,x_1,\textbf{u})\!-\!\textbf{x}(t,x_2,\textbf{v})|\leq\!\b_1(|x_1\!-\!x_2|,t)\!+\!\b_2(\|\textbf{u}\!-\!\textbf{v}\|)$ and parameters $\tau,\eta,\mu>0$
	satisfy the inequality
	\begin{equation}\label{eqnbisi}
	\beta_1(\e,\tau)+\beta_2(\mu)+\eta/2\leq \e,
	\end{equation}
	then $T_\tau(\S)\cong^{(\e,\mu)} T_{\tau\mu\eta}(\Sigma)$.
\end{lemma}

Lemma \ref{thmsi} can be proved following the proof of Theorem 5.1 in \cite{PolaAuto08}. The key step is to show (ii) and (iii) in Definition \ref{dfnemubi} by the following fact: for any input $\ell_1\in U_1$ of $T_\tau(\S)$, we can choose input $\ell_2\in U_2$ of $T_{\tau\mu\eta}(\Sigma)$ such that $|\ell_1-\ell_2|_\infty\leq \mu$, and vice versa. Specifically, we let $\ell_2=\mathcal{Q}(\ell_1)$ where the quantization function $\mathcal{Q}(\cdot)$ is defined entry-wisely as in \eqref{quanineq}.



We call the transition system $T_{\tau\mu\eta}(\S)$ a \emph{symbolic model}  for $\S$ where $T_\tau(\S)\cong^{(\e,\mu)} T_{\tau\mu\eta}(\S)$. One nice property of this symbolic model is that its evolution can be chosen to be deterministic \cite{girard2013low}, which is appropriate for discrete software-based implementation.
Particularly, if the state space and input space of $T_{\tau\m\eta}(\S)$ are chosen to be bounded sets, then the resulting $T_{\tau\m\eta}(\S)$ will be a finite transition system.

\subsection{Passivity Indices}\label{secdis}

\begin{definition}\label{ddissCT}\cite{BrogDissBook07}
	A continuous-time control system $\S$ as in \eqref{ctd} is called \emph{dissipative} with respect to a supply function $w(u,y)$ if there exists a positive semi-definite  storage function $V(x)$ such that the following (integral) dissipation inequality is satisfied for any  $x(t_0)\in\R^n$ with $t_0,t_1\in\R_{\geq 0},t_0<t_1$, and any admissible inputs $u$:
	\begin{equation}\label{con:dissipa}
	V(x(t_1))-V(x(t_0))\leq \int_{t_0}^{t_1}w(u(s),y(s))\,\ud s.
	\end{equation}
\end{definition}
If $V(x)$ is differentiable, an equivalent (differential) form of \eqref{con:dissipa} is $\dot{V}(x(t))\leq w(u(t),y(t)),\;\forall t\in\R_{\geq 0}.$

\begin{definition}\label{ddissDT}\cite{BrogDissBook07}
	A discrete-time control system $\S_d$ as in \eqref{dtd} is called \emph{dissipative} with respect to the supply function $w(u,y)$ if there exists a positive semi-definite storage function $V(x)$ such that  the following dissipation inequality is satisfied for any $x[k_0]\in\R^n$ with $k_0,k_1\in\Z_{\geq 0},k_0<k_1$,  and any admissible input $u$:
	\begin{equation}\label{dis:dissipa}
	V(x[k_1])-V(x[k_0])\leq \sum_{k=k_0}^{k_1-1}w(u[k],y[k]).
	\end{equation}
\end{definition}
{\color{black}The summation form of \eqref{dis:dissipa} is equivalent to the one-step form $V(x[k+1])-V(x[k])\leq w(u[k],y[k]),\;\forall k\in\Z_{\geq 0}$ \cite{BrogDissBook07}.
In the rest of the paper, the input and the output of $\S$ (or  $\S_d$) are assumed to have the same dimension (i.e., $m=s$) without loss of generality \cite{Khalil2000}.}

\begin{definition}\label{def:IFOFP}
	A system $\S$ (or $\S_d$) is called \emph{input feedforward output feedback passive} (IF-OFP) with passivity indices $(\nu,\rho)$ if it is dissipative with respect to the supply function $w(u,y)=u^\top y-\nu u^\top u-\rho y^\top y$
	for some $\nu,\rho\in\R$, denoted as IF-OFP$(\nu,\rho)$.
\end{definition}


Passivity indices $\nu,\rho$ reflect the excess or shortage of passivity of a system {\color{black}where positive values reflect excess and negative values reflect shortage. A passive system  is IF-OFP$(0,0)$, and a system that is 
IF-OFP$(\nu,\rho)$ with $\nu> 0,\rho> 0$ has an excess of passivity.}
The indices  $\nu,\rho$ are not unique because 
for a system that is IF-OFP$(\nu,\rho)$, it is also IF-OFP$(\nu',\rho')$ for any $\nu'<\nu,\rho'<\rho$.

\subsection{Problem Statement}
\begin{figure}[!hbt]
	\begin{subfigure}[b]{\linewidth}
		\centering
		\includegraphics[width=\linewidth]{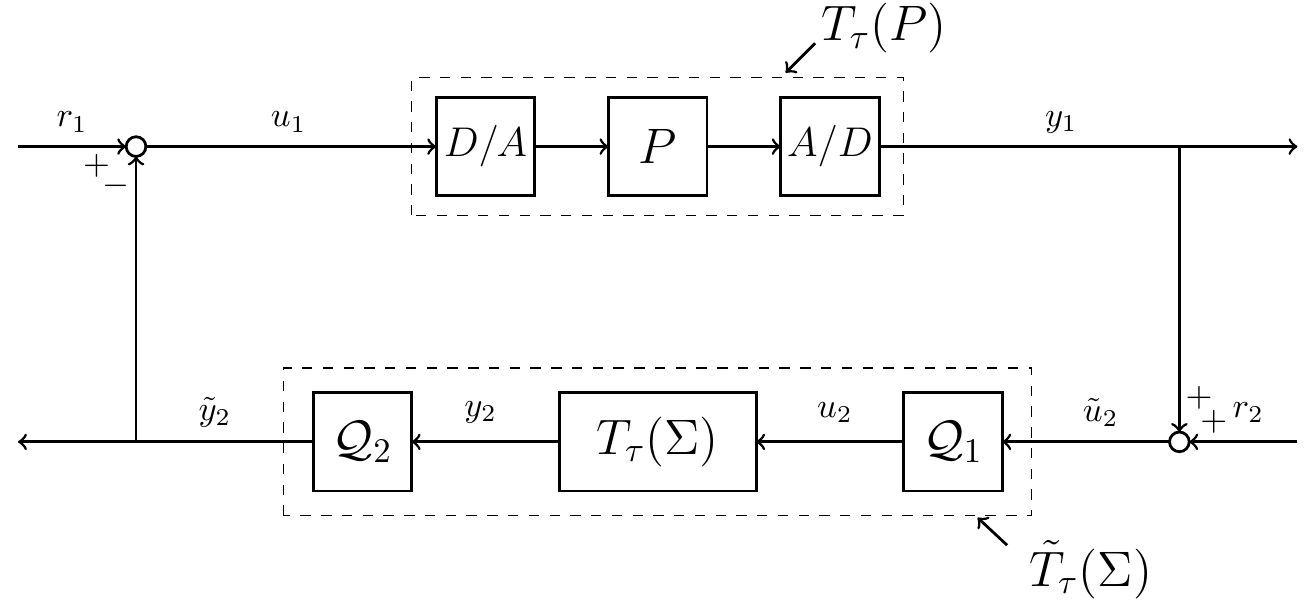}
		\caption{{\color{black} Feedback connection of $T_\tau(P)$ and $\tilde T_{\tau}(\S)$, denoted as  $T_\tau(P)\times_{\mathcal{F}}\tilde T_{\tau}(\S)$. }}
	\end{subfigure}\vskip 1mm
	\begin{subfigure}[b]{\linewidth}
		\centering
		\includegraphics[width=\linewidth]{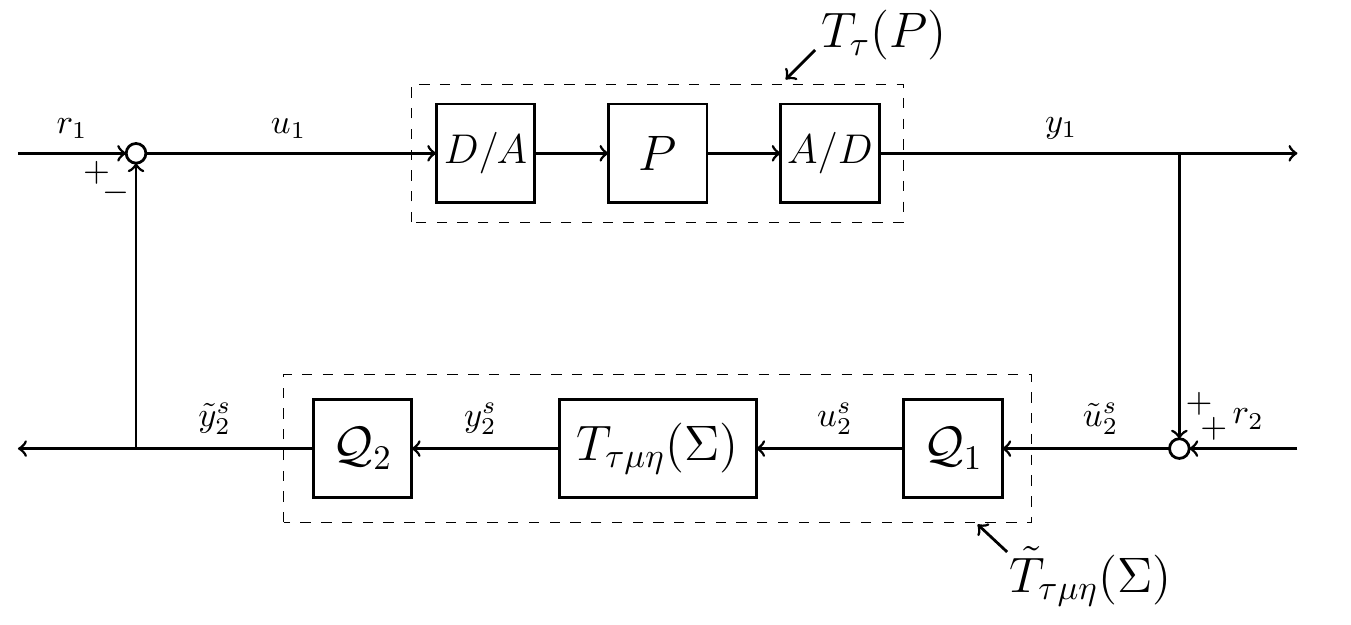}
		\caption{{\color{black} Feedback connection of $T_\tau(P)$ and $\tilde T_{\tau\mu\eta}(\S)$, denoted as $T_\tau(P)\times_{\mathcal{F}}\tilde T_{\tau\mu\eta}(\S)$.}}
	\end{subfigure}
	\caption{Two closed-loop system setups considered.}\label{figfe}
\end{figure}

We will consider the closed-loop configurations shown in Fig. \ref{figfe}, where the  feedback connection of $T_\tau(P)$ and $\tilde T_\tau(\S)$ is denoted as $T_\tau(P)\times_{\mathcal{F}}\tilde T_\tau(\S)$, and the  feedback connection of $T_\tau(P)$ and $\tilde T_{\tau\mu\eta}(\S)$ is denoted as $T_\tau(P)\times_{\mathcal{F}}\tilde T_{\tau\mu\eta}(\S)$. In these configurations, $P$ and $\S$ are both continuous models, which can be considered as the original plant and the controller, respectively; $T_\tau(P)$ (resp. $T_\tau(\S)$) is the time-sampled model of $P$ (resp. $\S$), which consists of $P$ (resp. $\S$), a zero-oder hold device and an uniform sampler; $\tilde T_\tau(\S)$ (resp. $\tilde T_{\tau\mu\eta}(\S)$) is the model that consists of $T_\tau(\Sigma)$ (resp. $T_{\tau\mu\eta}(\S)$), the input quantizer $\mathcal{Q}_1$ and the output quantizer $\mathcal{Q}_2$. The setup in Fig. \ref{figfe} (a) can be considered as replacing $\S$ with $\tilde T_\tau(\S)$, the sampled and quantized model of $\S$; the setup in Fig. \ref{figfe} (b) can be considered as replacing $\S$ with $\tilde T_{\tau\mu\eta}(\S)$, the approximate bisimilar symbolic model of $T_\tau(\S)$. We assume that 1) the external reference inputs $r_i, i=1,2$ to the closed-loop system are discrete-time signals;
2)  the sampling times of $T_\tau(P)$ and $T_\tau(\Sigma)$ are both $\tau$;  3) all the discrete-time signals in the feedback loop are synchronized.

{\color{black}
	The problems that will be studied are the following: 
	\emph{Given the passivity indices of $P$ and $\S$ and the setup of Fig. \ref{figfe} (a), what passivity property can be preserved for $T_\tau(P)\times_{\mathcal{F}}\tilde T_\tau(\S)$, and under what conditions the states of $T_\tau(P)\times_{\mathcal{F}}\tilde T_\tau(\S)$ are bounded? Similarly, given the passivity indices of $P$ and $\S$ and the setup of Fig. \ref{figfe} (b), under what conditions the states of $T_\tau(P)\times_{\mathcal{F}}\tilde T_{\tau\mu\eta}(\S)$ are bounded?} 
	
	
	We will study the passivity of $T_\tau(P)\times_{\mathcal{F}}\tilde T_\tau(\S)$ and give a boundedness result for it in Section \ref{sec:connect}, based on the passivity degradation results in Section \ref{sec:degradation}. After that, we will investigate the state boundedness of $T_\tau(P)\times_{\mathcal{F}}\tilde T_{\tau\mu\eta}(\S)$ in Section \ref{sec:symbolic}. The main theorems of the paper are summarized in Table \ref{tab1}.}

\begin{table}[!hbt]
	\begin{center}
		\begin{tabular}{c|c}
			\hline
			Theorem \ref{thm1}	& passivity degradation from $\Sigma$ to $T_\tau(\S)$ \\
			\hline 
			Theorem \ref{thm2}	& passivity degradation from $\Sigma_\tau$ to $\tilde T_\tau(\S)$  \\
			\hline
			Theorem \ref{thmfeedback}	& passivity inequality of $T_\tau(P)\times_{\mathcal{F}}\tilde T_{\tau}(\S)$  \\ 
			\hline
			Theorem \ref{thm3}	&  state boundedness of $T_\tau(P)\times_{\mathcal{F}}\tilde T_{\tau}(\S)$ \\ 
			\hline
			Theorem \ref{corsym}	&  state boundedness of $T_\tau(P)\times_{\mathcal{F}}\tilde T_{\tau\mu\eta}(\S)$ \\
			\hline
		\end{tabular}\caption{Summary of the main theorems. }\label{tab1}
	\end{center}
\end{table}
\section{Quasi-Passivity and  Strong Detectability}\label{sec:quasi}
In this section, we introduce quasi-passivity, strong detectability and some relevant lemmas, which are of interest by their own and will be  used in later sections.

\subsection{Input Feedforward Output Feedback Quasi-passivity}
{\color{black}Unlike dissipative systems, quasi-dissipative systems (see Def. 1 of \cite{PolQuasiTAC04}) or almost-dissipative  systems  (see Def. 2.1 of \cite{DowSCL03}) allow for ``internal energy generation''. The following definition is inspired by \cite{PolQuasiTAC04,DowSCL03} and Def. \ref{def:IFOFP}.}



\begin{definition}\label{defIFOFQP}
	The continuous (resp. discrete) time system $\S$ in \eqref{ctd} (resp. $\S_d$ in \eqref{dtd}) is called \emph{input feedforward output feedback quasi-passive} (IF-OFQP) if it is dissipative with respect to the supply function 
	$w(u,y)=u^\top y-\nu u^\top u-\rho y^\top y+\delta$, denoted as IF-OFQP$(\nu,\rho,\delta)$, {\color{black}where $\delta\geq 0$ is a constant.}
\end{definition}
{\color{black}By Def. \ref{defIFOFQP}, a discrete-time system $\S_d$ that is IF-OFQP$(\nu,\rho,\delta)$ satisfies the following  inequality}
\begin{align}
&V(x[k_1])-V(x[k_0])\nonumber\\
\leq& \sum_{i=k_0}^{k_1-1}\left( u[i]^\top y[i]-\nu |u[i]|^2-\rho |y[i]|^2+\delta\right)\label{ieq:VSQP}
\end{align}
for any $k_0,k_1\in\Z_{\geq 0},k_0<k_1$ and any admissible input $u[i]$ where $V(x)$ is a positive semi-definite function. {\color{black}Clearly, IF-OFP is a special case of IF-OFQP  with $\delta=0$.}

\begin{figure}[!h]
	\centering
	\includegraphics[width=0.6\linewidth]{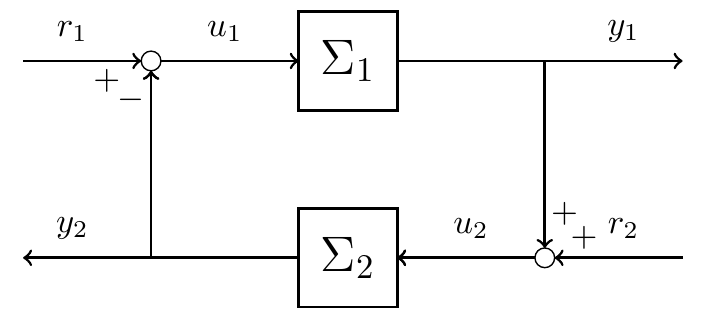}
	\caption{Feedback connection of two discrete-time systems}\label{figfeedback}
\end{figure}

{\color{black}
Consider the feedback-connection of two discrete-time systems $\S_1$, $\S_2$ shown in Fig. \ref{figfeedback}. Define $x_1$ and $x_2$ as the states, and $V_1(x_1)$ and $V_2(x_2)$ as the storage functions of $\S_1$ and $\S_2$, respectively. Define $r:=(r_1^\top,r_2^\top)^\top$ as the external input, and $y:=(y_1^\top,y_2^\top)^\top$ as the overall output. The following Lemma shows passivity indices of the feedback connection of $\S_1$ and $\S_2$. Its proof is similar to that of Theorem 6 in \cite{ZhuFeedbackACC14} and is thus omitted.}

\begin{lemma}\label{lemfeedind}
	Consider the feedback connection of two discrete-time systems $\S_1,\S_2$ as shown in Fig. \ref{figfeedback}. If $\S_i (i=1,2)$ is IF-OFQP$(\nu_i,\rho_i,\delta_i)$, then the feedback-connected system is IF-OFQP$(\hat \nu,\hat \rho,\hat\delta)$ with respect to the input $r$ and the output $y$  where $\hat{\delta}=\delta_1+\delta_2$ and $\hat\nu,\hat\rho$ can be chosen as
	\begin{equation}\label{feedineq1}
	\left\{
	\begin{aligned}
	&\hat\nu<\min\{\nu_1,\nu_2\},\\
	&\hat\rho\leq \min\{\rho_1-\frac{\hat \nu\nu_2}{\nu_2-\hat{\nu}},\rho_2-\frac{\hat \nu\nu_1}{\nu_1-\hat{\nu}}\}.
	\end{aligned}\right.
	\end{equation}
\end{lemma}

\subsection{Strong Detectability}
\begin{definition}\label{dfndete}
	{\color{black}The discrete-time system (\ref{dtd}) is said to be \emph{$N$-step strongly detectable} (SD) where $N\in\Z_{\geq 0}$, if there exist a constant $\vartheta\geq 0$ and a radially unbounded, positive definite function $p: \R^n\rightarrow \R_{\geq 0}$ such that for any $k_0\in\Z_{\geq 0}$, any initial state $x[k_0]\in\R^n$ and any admissible input $u[k]$, the following holds:
	\begin{equation}
\sum_{k=k_0}^{k_0+N}\vartheta|u[k]|^2+|y[k]|^2\geq p(x[k_0]).\label{strongdetec}
	\end{equation}}
\end{definition}

The intuition behind strong  detectability is large initial states must dictate large input-output signals. {\color{black}Similar definitions were also given in Eq. (38) in \cite{sontag1989smooth} and Def. 3 in \cite{PolQuasiTAC04}.} 
The strong  detectability in Def. \ref{dfndete} implies the zero-state observability given in Def. 6.5 of \cite{Khalil2000}, because  $u[k]={\bf 0},y[k]={\bf 0},k=k_0,...,k_0+N_0$, implies $x[k_0]={\bf 0}$. 

\begin{remark}
	A discrete-time linear system that is observable with an observability index $v$ is $v$-step SD. Given a system $x[k+1]=A_dx[k]+B_du[k],y[k]=C_dx[k]+D_du[k]$, define $U=[u[k_0]^\top,...,u[k_0+v]^\top]^\top$, $Y=[y[k_0]^\top,...,y[k_0+v]^\top]^\top$. The matrix $Y$ can be expressed as $Y=Ox[k_0]+HU$ where  $O$ is the observability matrix and $H$ is a matrix that can be constructed from $A_d,B_d,C_d,D_d$. Hence,  $\sum_{k=k_0}^{k_0+v}\vartheta|u[k]|^2+|y[k]|^2=\vartheta U^\top U+x[k_0]^\top O^\top Ox[k_0]+U^\top H^\top HU+2x[k_0]^\top O^\top HU$. Since $rank(O^\top O)=rank(O)=n$, $O^\top O\succ 0$. Then it is easy to see that there exist $\vartheta>0$ and {\color{black}$p(x)=x^\top Mx$ where $M\succ0$} such that \eqref{strongdetec} holds. 
\end{remark}

\begin{example}
	Consider the double integrator system $x_1[k+1]=x_2[k]$, $x_2[k]=u[k]$, $y[k]=x_1[k]+u[k]$ where $x_1,x_2,y,u\in \R$. It is easy to verify that the system is not $0$-step SD, but  $1$-step SD with ${\color{black}\vartheta}=2$ and {\color{black}$p(x)=\frac{1}{2}(x_1^2+x_2^2)$}.
\end{example}

The strong detectability of a feedback-connected system can be derived by the strong detectability of each subsystem, which is shown by the following Lemma \ref{lemstrongdet}. The proof of this lemma is given in Appendix \ref{sub:strongdetec}.
{\color{black}
\begin{lemma}\label{lemstrongdet}
Consider the feedback connection of two discrete-time systems $\S_1,\S_2$ as shown in Fig. \ref{figfeedback}. Suppose that $\S_1$ is $N_1$-step SD that satisfies  
\begin{align}
\sum_{k=k_0}^{k_0\!+\!N_1}\vartheta_1|u_1[k]|^2\!+\!|y_1[k]|^2\!\geq\! p_1(x_1[k_0]),\forall k_0\in\Z_{\geq 0},\label{det1}
\end{align}
and $\S_2$ is $N_2$-step SD that satisfies 
\begin{align}
\sum_{k=k_0}^{k_0\!+\!N_2}\vartheta_2|u_2[k]|^2\!+\!|y_2[k]|^2\!\geq\! p_2(x_2[k_0]),\forall k_0\in\Z_{\geq 0},\label{det2}
\end{align}
where $\vartheta_1,\vartheta_2\geq 0$ and $p_1,p_2$ are radially unbounded, positive definite functions, then, the feedback-connected system is $N$-step SD with respect to the input $r$ and the output $y$ where $N=\max\{N_1,N_2\}$, $p(x)=(1-\vartheta)(p_1(x_1)+p_2(x_2))$ and
\begin{align}\label{vartheta}
\vartheta=\max\{\frac{2\vartheta_1}{2\vartheta_1+1},\frac{2\vartheta_2}{2\vartheta_2+1}\}.
\end{align} 
\end{lemma}}

The next lemma shows that a discrete-time system $\S_d$ that is SD and IF-OFQP with $\rho>0$ has the (uniform) bounded-input-bounded-state property. The proof of this lemma is given in Appendix \ref{sec:prf0}.

\begin{lemma}\label{lembound}
	Suppose that the system  $\S_d$ given in \eqref{dtd} is 1) $N$-step SD satisfying \eqref{strongdetec}, 2) IF-OFQP$(\nu,\rho,\delta)$ satisfying \eqref{ieq:VSQP} with $\rho>0,\delta\geq 0$ and a function $V(x)$ that is continuous, positive semi-definite, radially unbounded. Let $\l$ be a number such that $0<\l<\rho$, and define $\eta_1=\frac{1}{4\l}-\nu> 0$, $\eta_2=\rho-\l>0$. Then, \\
	i) for any $k\in\Z_{\geq 0}$, it holds that $x[k]\in\mathcal{D}_1$ where 
	\begin{align}
	{\color{black}\mathcal{D}_1:=\{z\mid V(z)\leq \xi_1\}}
	\end{align}
	and ${\color{black}\xi_1=c_1+c_3}$, ${\color{black}\xi_2=\max\{p(x[0]),c_2/\eta_2\}}$, $c_1=\max_{z\in\mathcal{C}_1}V(z)$, $c_2=(N+1)[(\eta_1+\vartheta\eta_2)\|u\|^2+\delta]$, $c_3=(N+1)(\eta_1\|u\|^2+\delta)$, and ${\color{black}\mathcal{C}_1=\{z\mid p(z)\leq \xi_2\}}$.

	ii) there exists $K\in\Z_{\geq 0}$ such that $x[k]\in\mathcal{D}_2$ for all $k\geq K$ where  
	\begin{align}
	{\color{black}\mathcal{D}_2:=\{x\mid V(x)\leq \xi_3\}}\label{D2}
	\end{align}
	and $\xi_3=c_2+c_4$, $\xi_4=(c_2+ c_5)/\eta_2$, $c_2$ is as in 1), $c_4=\max_{z\in\mathcal{C}_2}V(z)$, $c_5$ is a postive number, and ${\color{black}\mathcal{C}_2=\{z\mid p(z)\leq \xi_4\}}$.

\end{lemma}

\begin{remark}
	Lemma 6.7 in \cite{Khalil2000} shows the zero-input Lyapunov stability for continuous-time, output  strictly passive systems with the zero-state  observability assumption. {\color{black}Lemma \ref{lembound} is a complement to Lemma 6.7 in the discrete-time, quasi-passivity setting.}
\end{remark}




\subsection{Input-to-state Practical Stability}
\begin{definition}\cite{lazar2008input}\label{dfnISPS}
	The discrete-time system $\S_d$ given in \eqref{dtd} is called (globally) \emph{input-to-state practically stable} (ISpS) if there exist a $\mathcal{KL}$-function $\beta_1$, a $\mathcal{K}$-function $\beta_2$ and a positive constant $d$ such that, for any $u[k]$ with $\|u\|<\infty$ and any $x_0\in\R^n$, it holds that 
	$$
	|x[k]|\leq \beta_1(|x_0|,k)+\beta_2(\max_{0\leq j\leq k-1}|u[j]|)+d,\;\forall k\in\Z_{\geq 0}.
	$$ 
\end{definition}

\begin{definition}\cite{lazar2008input}
	A continuous function $V:\R^n\rightarrow\R_{\geq 0}$ is called the \emph{ISpS-Lyapunov function} for $\S_d$ if there exist $d_1,d_2,a,b,c,\l>0$ with $c\leq b$ and a
	$\mathcal{K}$-function $\sigma$, such that
	the following hold: 1) 
	$\alpha_1(|x|)\leq V(x)\leq \alpha_2(|x|)+d_1,\forall x\in\R^n$ where $\alpha_1(s)=as^\l$, $\alpha_2(s)=bs^\l$; 2) 
	$V(f_d(x,u))-V(x)\leq -\alpha_3(|x|)+\sigma(|u|)+d_2$, $\forall x\in\R^n,\forall u\in\R^m$,
	where $\alpha_3(s)=cs^\l$.
\end{definition}
Clearly, $\S_d$ is input-to-state stable (ISS) when $d=0$ in Def. \ref{dfnISPS}. 
In Theorem 2.5 in \cite{lazar2008input}, it was shown  that the discrete-time ISpS-Lyapunov function is a sufficient condition for $\S_d$ to be ISpS. 

The following Lemma \ref{lem4} shows conditions under which a discrete-time, quasi-passive system $\S_d$ is ISpS. Its proof is shown in Appendix \ref{prooflem4}.

\begin{lemma}\label{lem4}
	If there exist $\lambda,a,b,c>0,d_1\geq 0$ such that the system  $\S_d$ given in \eqref{dtd} is 1) $0$-step SD satisfying \eqref{strongdetec} with {\color{black}$p(x)\geq c|x|^\lambda$}, 2) IF-OFQP$(\nu,\rho,\delta)$ satisfying \eqref{ieq:VSQP} with $\rho>0,\delta\geq 0$ and a continuous function $V(x)$ that satisfies $a|x|^\lambda\leq V(x)\leq b|x|^\lambda+d_1$,  then $\S_d$ is ISpS.
\end{lemma}

Lemma \ref{lem4} establishes the connection between quasi-passivity and ISpS, which enables us to use the results of ISpS to analyze IF-OFQP systems under certain circumstances \cite{JiangAuto96,JiangISSAuto01,lazar2008input}.
For instance, if $V(x)$ in Lemma \ref{lembound} satisfies $\alpha_1(|x|)\leq V(x)\leq \alpha_2(|x|),\forall x\in\R^n$, where $\alpha_1,\alpha_2$ are  $\mathcal{K}_\infty$-functions, then the ultimate bound of $\S_d$ can be derived similar to the proof of Lemma 3.5 in \cite{JiangISSAuto01}. Specifically, the set $\mathcal{D}_2$ can be given as 
\begin{align}\label{boundD2}
{\color{black}\mathcal{D}_2=\{x\mid V(x)\leq \xi_5\}}
\end{align}
{\color{black}where $\xi_5=\alpha_4^{-1}\circ\alpha_5^{-1}\big((\eta_1+\vartheta\eta_2)\|u\|^2+\delta\big)$,}
$\alpha_4\in\mathcal{K}_\infty$ satisfies $\alpha_4\leq \alpha_3\circ\alpha_2^{-1}$,  $id-\alpha_4\in\mathcal{K}$, $id$ is the identity function, $\alpha_3(x)=p(x)\in\mathcal{K}_\infty$, and $\alpha_5\in\mathcal{K}_\infty$ satisfies $id-\alpha_5\in\mathcal{K}_\infty$. 
{\color{black}Note that the set $\mathcal{D}_2$ expressed in \eqref{D2} or \eqref{boundD2} 
	can be made arbitrarily small by choosing $\|u\|$ and $\delta$ small enough.}

\begin{remark}
	In \cite{gross2016relaxed}, the input-to-state stability was discussed by a relaxed ISS-Lyapunov function, which satisfies a weaker decrease condition defined over a bounded time interval. Because of the N-step strongly detectability condition, the proof of ultimate boundedness in  Lemma \ref{lembound} also has to consider the decrease of $V$ every N+1 steps; however, it is different from the proof of \cite{gross2016relaxed}.
\end{remark}


\begin{figure}[!b]
	\begin{subfigure}[b]{\linewidth}
		\centering
		\includegraphics[width=0.3\linewidth]{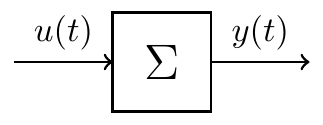}
		\caption{$\Sigma$ is a continuous time model. }
	\end{subfigure}\vskip 1mm
	\begin{subfigure}[b]{\linewidth}
		\centering
		\includegraphics[width=0.7\linewidth]{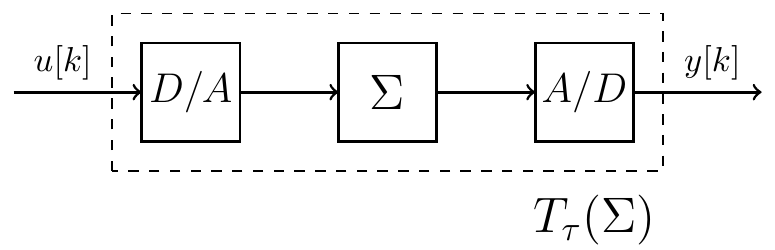}	
		\caption{$T_\tau(\Sigma)$ is a time-sampled model of $\Sigma$, which consists of $\Sigma$, a ZOH device and an uniform sampler.}
	\end{subfigure}\vskip 1mm
	\begin{subfigure}[b]{\linewidth}
		\centering
		\includegraphics[width=0.7\linewidth]{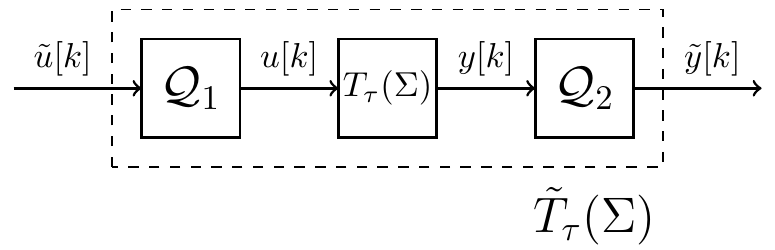}
		\caption{$\tilde T_\tau(\Sigma)$ is the model that consists of $T_\tau(\Sigma)$, the input quantizer $\mathcal{Q}_1$ and the output quantizer $\mathcal{Q}_2$.}
	\end{subfigure}
	\caption{Three system setups considered.}\label{fig1}
\end{figure}

\section{Passivity Degradation Under Sampling and Quantization}\label{sec:degradation}

In this section, we study the degradation of passivity indices from a continuous-time system $\Sigma$ to the time-sampled system $T_\tau(\S)$ and the time-sampled, input/output quantized system $\tilde T_{\tau}(\Sigma)$. Configurations of $\Sigma$, $T_\tau(\Sigma)$ and $\tilde T_\tau(\Sigma)$ are shown in Fig. \ref{fig1}.

We make two assumptions on the output function of $\S$.
\begin{assumption}\label{assadd}
	The output function $h$ of $\S$ has the additive form $h(x,u)=h_1(x)+h_2(u)$.
\end{assumption}
\begin{assumption}\label{assgain}
	There exist a constant $\gamma>0$ and a function $\beta:\R^n\to \R_{\geq 0}$ with $\beta({\bf 0})=0$ such that for any initial state $x_0\in\R^n$, $h_1$ satisfies the following inequality for any $T>0$ and any admissible $u(t)$:
	\begin{equation}\label{eqngain}
	\int_{0}^{T}|\dot{h}_1(x(t))|^2\,\ud t\leq\gamma^2 \int_{0}^{T}|u(t)|^2\,\ud t+\beta(x_0).
	\end{equation}
\end{assumption}
\begin{remark}
{\color{black}	A gain assumption without considering the initial state in \eqref{eqngain} was made in Theorem 2 of \cite{OishiCDC10} and Assumption 2 of \cite{xia2017approx}}: for any $T>0$ and any admissible $u(t)$, it holds that 
	\begin{align}
	\int_{0}^{T}|\dot{h}(t)|^2\,\ud t\leq\gamma^2 \int_{0}^{T}|u(t)|^2\,\ud t.\label{gainassum}
	\end{align}
	But for the inequality \eqref{gainassum} to hold, the function $h$ in \cite{OishiCDC10,xia2017approx} needs to be independent of $u$ and additionally, the initial state is assumed to be zero (i.e., $x_0={\bf 0}$). In comparison, Assumption \ref{assadd} and Assumption \ref{assgain} are less restrictive by considering the additive form  and  the initial condition explicitly.
\end{remark}
\begin{remark}\label{remassgain}
	If there exist a positive semi-definite function $\beta(x)$ and a positive number $\gamma$ such that 
	\begin{align}\label{gainlya}
	\dot{\beta}(x(t))\leq \gamma^2|u(t)|^2-|\dot{h}_1(t)|^2,\;\forall t\in \R_{\geq 0},
	\end{align}
	then Assumption \ref{assgain} holds. Particularly,
	for a stable LTI system $\dot{x}=Ax+Bu$ with $h_1(x)=Cx$,
	if there exist a matrix $P\succeq 0$ and a number $\gamma>0$ such that 
	\begin{equation*}
	\left(\begin{array}{cc}
	A^\top P+PA+A^\top C^\top CA &\;\; PB+A^\top C^\top CB \\ 
	B^\top P+B^\top C^\top CA &\;\;-\gamma^2 I+B^\top C^\top CB
	\end{array} \right)\preceq 0
	\end{equation*}
	{\color{black}then one can verify that the function $\beta(x)=x^TPx$ satisfies \eqref{gainlya}, implying that Assumption \ref{assgain} holds for such a system. }
	For nonlinear systems whose dynamics are polynomial,  the sum-of-squares optimization can be used to search for $\beta(x)$ that satisfies \eqref{gainlya} \cite{xia2017approx}. 
\end{remark}

\subsection{Passivity Degradation of $T_{\tau}(\Sigma)$}
In this subsection, we consider deriving passivity indices of $T_{\tau}(\S)$ from those of $\S$. As shown in Fig. \ref{fig1} (a), $\S$ is a continuous-time system whose input and output are $u(t)$ and $y(t)$, respectively, and $T_{\tau}(\Sigma)$  is a discrete-time system whose input and output are $u[k]$ and $y[k]$, respectively.
The following theorem shows quantitatively how passivity indices degrade from $\S$ to $T_\tau(\S)$. The proof of the theorem is given in Appendix \ref{sec:prf2}.

\begin{theorem}\label{thm1}
	Suppose that $\S$ satisfies Assumption \ref{assadd} and Assumption \ref{assgain}, and is IF-OFP$(\nu,\rho)$ with a positive semi-definite storage function $V(x)$. Then $T_\tau(\S)$ satisfies the following inequality  for any $k_0,k_1\in\Z_{\geq 0},k_0<k_1$ and admissible inputs $u\in U_1$:
	\begin{equation}\label{thm1:ineq1}
	\begin{aligned}
	\hat V(x[k_1])-\hat V(x[k_0]) &\leq \delta(x[k_0])+\sum_{k=k_0}^{k_1-1} \left(u[k]^\top y[k]\right.\\
	&\quad \left.-\nu'|u[k]|^2-\rho'|y[k]|^2\right)
	\end{aligned}
	\end{equation}
	where $\hat V(x)=\frac{1}{\tau}V(x)$ and
	\begin{equation}\label{indextau1}
	\left\{
	\begin{aligned}
	&\nu'=\nu - \tau\gamma-\tau^2\gamma^2(1+\lambda_1)|\rho|,\\
	&\rho'=\rho - |\rho|/\lambda_1,\\
	&\delta(x[k_0])=w\beta(x[k_0]),\\
	&w=|\rho|\tau(1+\l_1)+\frac{1}{\gamma},
	\end{aligned}\right.
	\end{equation}
	with $\l_1>0$ an arbitrary positive number.
\end{theorem}




Although inequality \eqref{thm1:ineq1} has an additional bias term $\delta$ on its right-hand side, we still call \eqref{thm1:ineq1} the \emph{passivity inequality} satisfied by $\S_d$, with some abuse of language.
Note that this bias term $\delta$ will not affect the summation of the supply function, and $\delta=0$ when $x[0]={\bf 0}$.

\begin{remark}
Equation \eqref{indextau1}  indicates that $\nu'<\nu$, $\rho'<\rho$, and when $\lambda_1$ is fixed, a smaller sampling time $\tau$ can result in a smaller passivity degradation (i.e., larger $\nu',\rho'$), which was also shown in \cite{OishiCDC10,xia2017approx}. Theorem \ref{thm1} also generalizes the corresponding results of \cite{OishiCDC10,xia2017approx} in several ways.
Firstly, a more general gain assumption \eqref{eqngain} that takes into account the initial condition is used, which induces  a bias term  $\delta(x[0])$ on the right hand side of \eqref{thm1:ineq1}. Secondly, the positive number $\l_1$ in \eqref{indextau1}  provides additional flexibility in balancing the passivity indices $\nu',\rho'$ and the term $\delta$.
Furthermore, to obtain the degraded indices in \cite{xia2017approx}, it was assumed   that $\nu>0$ and an inequality about $\tau,\nu,\rho$ needs to be satisfied (see Corollary 6 in \cite{xia2017approx}), but neither of them is required in Theorem \ref{thm1}.
\end{remark}

\begin{remark}\label{remark4}
	The system $\S$ in Theorem \ref{thm1} can be nonlinear control system in general. 
	When $\S$ is a continuous-time LTI system, dynamics of $T_\tau(\S)$ can be given explicitly and 
	passivity indices of $T_\tau(\S)$ can be obtained directly by solving an LMI {\color{black}(e.g., see Lemma 3 in \cite{kottenstette2014relationships}). Note that in this case the term $\delta$ does not exist}.
\end{remark}

\begin{example}\label{ex1}
	Consider the following LTI system:
	\begin{align*}
	\S:
	\begin{cases}\dot{x}=\left[\begin{array}{cc}-1.8&-1.3\\
	1.2&-2.5
	\end{array}\right]x+
	\left[\begin{array}{cc}0.2& 0\\0&0.3
	\end{array}\right]u,\\
	y=\left[\begin{array}{cc}0.2& -0.3\\0.3&0.15
	\end{array}\right]x + \left[\begin{array}{cc}0.5& 0\\0&0.4
	\end{array}\right]u,
	\end{cases}
	\end{align*}
	where $x(t),u(t),y(t)\in\R^2$. 
	One can verify that $\S$ is IF-OFP$(0.3,0.5628)$ with a positive definite storage function 
	$V(x)= 0.23|x|^2$.
	Furthermore, Assumption \ref{assadd} holds and $h_1$ satisfies Assumption \ref{assgain} with $\gamma=0.2$ and $\beta(x)= 0.2187|x|^2$.
	{\color{black}Suppose that the sampling time $\tau=0.3$. 	
	Letting $\l_1=10$ in \eqref{indextau1}, we obtain $\nu'=0.2177,\rho'=0.5065$ and $\delta(x)=6.8572 \beta(x)$. By Theorem \ref{thm1}, $T_\tau(\S)$ satisfies 
	\begin{align*}
	&\hat V(x[k_1])-\hat V(x[k_0]) \leq 1.4997|x[k_0]|^2\\
	&\quad+\sum_{k=k_0}^{k_1-1} \left(u[k]^\top y[k]\!-\!0.2177|u[k]|^2 \!-\! 0.5065|y[k]|^2\right)
	\end{align*}
	for any $k_0,k_1\in\Z_{\geq 0},k_0<k_1$ 
	where $\hat V(x)=\frac{1}{\tau}V(x)=0.7667|x|^2$. By solving the LMI in Lemma 3 of \cite{kottenstette2014relationships}, one can verify that $T_\tau(\S)$ is IF-OFP$(0.20,0.9803)$ with a positive definite storage function $V(x)=0.23|x|^2$, which is deliberately chosen the same as above. Then, $T_\tau(\S)$ satisfies the following inequality 	for any $k\in\Z_{\geq 0}$: 
	\begin{align*}
	\hat V(x[k+1]) - \hat V(x[k]) &\leq u[k]^\top y[k] -0.2|u[k]|^2\\
	&\quad - 0.9803|y[k]|^2
	\end{align*}
	where $\hat V(x)=0.7667|x|^2$.
	From the two passivity inequalities above, one can observe that passivity indices derived from \eqref{indextau1} and from solving the LMI are not comparable directly.}

\end{example}

\subsection{Passivity Degradation of $\tilde T_{\tau}(\S)$}\label{subdegradtidT}
In this subsection, we consider deriving passivity indices of $\tilde T_{\tau}(\S)$ from those of $T_{\tau}(\S)$. Recall that the passivity indices of $T_\tau(\S)$ not only can be obtained from a continuous-time system $\S$ by Theorem \ref{thm1}, in which case $T_\tau(\S)$ satisfies \eqref{thm1:ineq1} with $\nu',\rho',\delta$ given in \eqref{indextau1}, but also may be obtained directly from the dynamics of  $T_\tau(\S)$, in which case the bias term $\delta$ in \eqref{thm1:ineq1} does not exist.

Note that the system $\tilde T_{\tau}(\S)$ shown in Fig. \ref{fig1} consists of $T_{\tau}(\S)$ and two uniform quantizers $\mathcal{Q}_1,\mathcal{Q}_2$. 
The (uniform) quantization function $\mathcal{Q}(\cdot)$ with a quantization precision $\mu$  is defined as 
\begin{equation}\label{quanineq}
\mathcal{Q}(s):=\left\{
\begin{aligned}
&\left\lfloor \frac{s}{\mu} \right\rfloor \mu, \;\;\mbox{if}\;s\geq 0;\\
&\left\lceil \frac{s}{\mu} \right\rceil \mu, \;\;\mbox{if}\;s< 0,
\end{aligned}
\right.
\end{equation}
where $\lfloor\cdot \rfloor$ is the floor function (i.e.,  $\lfloor x \rfloor$ is the greatest integer no larger than $x$) and $\lceil \cdot \rceil$ is the ceiling function (i.e., $\lceil x \rceil$ is the least integer no less than $x$). If $s$ is a vector, then $\mathcal{Q}$ is implemented entry-wisely.

The input of $\tilde T_{\tau}(\S)$, $\tilde u[k]$, and the quantized input, $u[k]$, are related by $u[k]=\mathcal{Q}_1(\tilde u[k])$ where the function $\mathcal{Q}_1(\cdot)$ is given in \eqref{quanineq} with a quantization precision $\mu_1$. It is easy to see that for any $k \in\Z_{\geq 0}$, $|u[k]| \leq  |\tilde u[k]|$ and $|u[k]-\tilde u[k]|\leq \sqrt{m}\mu_1$, which results from the facts that $|u[k]-\tilde u[k]|\leq \sqrt{m}|u[k]-\tilde u[k]|_\infty$ and $|u[k]-\tilde u[k]|_\infty\leq \mu_1$.
The output of $\tilde T_{\tau}(\S)$, $\tilde y[k]$, and the output before quantization, $y[k]$, are related by $\tilde y[k]=\mathcal{Q}_2(y[k])$ where the function $\mathcal{Q}_2(\cdot)$ is given in \eqref{quanineq} with a precision $\mu_2$. It is clear that $|y[k]-\tilde y[k]|\leq \sqrt{m}\mu_2$  for any $k \in\Z_{\geq 0}$.




The following theorem shows how passivity indices degrade from  $T_\tau(\S)$ to $\tilde T_{\tau}(\S)$. The proof of the theorem is given in Appendix  \ref{sec:pfdegrade}.


\begin{theorem}\label{thm2}
	Suppose that $T_\tau(\S)$  is IF-OFP$(\nu',\rho')$ that satisfies \eqref{thm1:ineq1} with a positive semi-definite storage function $\hat V(x)$. 
	Then $\tilde T_\tau(\S)$ satisfies the following passivity inequality  for any $k_0,k_1\in\Z_{\geq 0},k_0<k_1$ and any admissible inputs $u\in U_1$:
	\begin{align}
	&\hat V(x[k_1])-\hat V(x[k_0]) \leq \delta(x[k_0])+\sum_{k=k_0}^{k_1-1} \Big(\tilde u[k]^\top \tilde y[k]\big.\nonumber\\
	&\quad\quad\big.-\tilde\nu|\tilde u[k]|^2-\tilde\rho|\tilde y[k]|^2+\tilde \delta\Big)\label{thm2:ineq1}
	\end{align}
	where {\color{black}$\delta(x[k_0])$ is given in \eqref{indextau1},}
	\begin{equation}\label{indextau3}
	\left\{
	\begin{aligned}
	\tilde \nu&=\nu'-|\nu'|/\l_3-1/4\l_4,\\
	\tilde \rho&=\rho'-|\rho'|/\l_2-1/4\l_5,\\
	\tilde \delta&=[|\rho'|(1+\l_2)+\l_4]m\mu_2^2\!+\![|\nu'|(1+\l_3)+\l_5]m\mu_1^2,
	\end{aligned}\right.
	\end{equation}
	and  $\l_2,\cdots,\l_5>0$ are  arbitrary positive  numbers.
\end{theorem}

\begin{remark}\label{remark8}
{\color{black}In \eqref{thm2:ineq1}, the term $\tilde \delta$ exists 
	due to the input and output quantizations of $\tilde T_{\tau}(\S)$, and the term  $\delta$ is equal to $0$ if $\delta$ shown in \eqref{thm1:ineq1} is equal to $0$ (e.g., $\delta=0$ when $\nu,\rho$ are solved by an LMI in Example \ref{ex1}).} 
The following observations can also be made from \eqref{thm2:ineq1}.
Firstly,  when $\mu_1$ and $\mu_2$, the quantization precisions, are fixed, there is a trade-off in choosing $\l_i (2\leq i\leq 5)$: larger $\l_i$  result in larger $\tilde\nu$ and $\tilde\rho$ (i.e., more ``excess of passivity''), while smaller $\l_i$ result in smaller $\tilde \delta$ (smaller ``internal energy generation''). Secondly, if $\nu'>0$ (resp. $\rho'>0$), then it is always possible to choose $\l_2,\l_3$ (resp. $\l_4,\l_5$) large enough such that $\tilde \nu>0$ (resp. $\tilde \rho>0$). Thirdly, with fixed $\l_i$, smaller $\mu_1,\mu_2$ (i.e., more precise quantizations) result in smaller $\tilde \delta$ (i.e., less quantization effect), and $\tilde \delta\rightarrow 0$ when  $\mu_1,\mu_2\rightarrow 0$. Therefore, if $\nu'>0,\rho'>0$, then it is always possible to choose $\l_i (2\leq i\leq 5)$ large enough and $\mu_1,\mu_2$ small enough such that $\tilde \nu>0$, $\tilde \rho>0$ and $\tilde \delta$ arbitrarily small.
\end{remark}

\begin{example}\label{ex2}
	Consider the system  $T_\tau(\Sigma)$ in Example \ref{ex1} again, which is IF-OFP$(0.20,0.9803)$ with a storage function $\hat V(x)=0.7667|x|^2$ when  $\tau=0.3$. Suppose that $\mu_1=\mu_2=0.01$, and let $\l_i=20$ for $i=2,3,4,5$  in \eqref{indextau3}. Then, by equation \eqref{indextau3}, $\tilde T_{\tau}(\S)$ is IF-OFQP$(0.1775,0.9188,0.0130)$. 
	That is, $\hat V(x[k+1]) - \hat V(x[k]) \leq \tilde u_2[k] \tilde y[k] -0.1775|\tilde u[k]|^2 - 0.9188|\tilde y[k]|^2 + 0.0130$ for any $k\in\Z_{\geq 0}$.
\end{example}

\section{Passivity Analysis of the Closed-Loop System}\label{sec:connect}



\subsection{Passivity Property of $T_\tau(P)\times_{\mathcal{F}}\tilde T_\tau(\S)$}\label{subpassineq}

In this subsection, we analyze the passivity of the system $T_\tau(P)\times_{\mathcal{F}}\tilde T_\tau(\S)$ shown in Fig. \ref{figfe} (a). 

Suppose that $P$ satisfies Assumption \ref{assadd} and Assumption \ref{assgain} with constant $\gamma_1$ and function $\beta_1$, and $P$ is  IF-OFP$(\nu_1,\rho_1)$ with a continuous, positive definite, radially unbounded storage function  $V_1$. Denote the state, input  and output of $T_{\tau}(P)$ by $x_1[k]\in\R^{n_1},u_1[k]\in\R^{m},y_1[k]\in\R^{m}$, respectively. By Theorem \ref{thm1}, $T_{\tau}(P)$ satisfies the following inequality for any $k_0,k_1\in\Z_{\geq 0},k_0<k_1$:
\begin{align}
&\frac{1}{\tau}V_1(x_1[k_1])-\frac{1}{\tau} V_1(x_1[k_0]) \leq \delta_1(x_1[k_0])\nonumber\\
&\quad +\sum_{k=k_0}^{k_1-1} \big(u_1[k]^\top y_1[k]-\nu_1'|u_1[k]|^2-\rho_1'|y_1[k]|^2\big) \label{con:ineq1}
\end{align}
where
\begin{equation}\label{comind1}
\left\{
\begin{aligned}
&\nu_1'=\nu_1 - \tau\gamma_1-\tau^2\gamma_1^2(1+\lambda_{11})|\rho_1|,\\
&\rho_1'=\rho_1 - |\rho_1|/\lambda_{11},\\
&\delta_1(x_1[k_0])=w_1\beta_1(x_1[k_0]),\\
&w_1=|\rho_1|\tau(1+\lambda_{11})+1/\gamma_1,
\end{aligned}
\right.
\end{equation}
and $\lambda_{11}>0$ is a positive number.


Suppose that $\Sigma$ satisfies Assumption \ref{assadd} and Assumption \ref{assgain} with constant $\gamma_2$ and function $\beta_2$. Suppose also that $\S$ is  IF-OFP$(\nu_2,\rho_2)$ with a continuous, positive definite, radially unbounded storage function  $V_2$. Denote the state, input and output of {\color{black}$\tilde T_{\tau}(\Sigma)$} by $x_2[k]\in\R^{n_2},\tilde u_2[k]\in\R^{m},\tilde y_2[k]\in\R^{m}$, respectively. By {\color{black}Theorem \ref{thm1} and Theorem \ref{thm2}},  $\tilde T_{\tau}(\S)$ satisfies the following  inequality for any $k_0,k_1\in\Z_{\geq 0},k_0<k_1$:
\begin{align}
&\frac{1}{\tau}V_2(x_2[k_1])-\frac{1}{\tau} V_2(x_2[k_0]) \leq \delta_2(x_2[k_0])+ \nonumber\\
&\sum_{k=k_0}^{k_1-1}\Big(\tilde u_2[k]^\top \tilde y_2[k]-\tilde \nu_2|\tilde u_2[k]|^2-\tilde\rho_2|\tilde y_2[k]|^2+\tilde \delta_2\Big) \label{con:ineq2}
\end{align}
where
\begin{equation}\label{comind2}
\left\{
\begin{aligned}
&\tilde \nu_2=\nu_2'-|\nu_2'|/\lambda_{23}-1/4\lambda_{24},\\
&\tilde \rho_2=\rho_2'-|\rho_2'|/\lambda_{22}-1/4\lambda_{25},\\
&\tilde \delta_2=[|\rho_2'|(1+\lambda_{22})+\lambda_{24}]m\mu_2^2\!+\![|\nu_2'|(1+\lambda_{23})+\lambda_{25}]m\mu_1^2,\\
&\nu_2'=\nu_2 - \tau\gamma_2-\tau^2\gamma_2^2(1+\lambda_{21})|\rho_2|,\\
&\rho_2'=\rho_2 - |\rho_2|/\lambda_{21},\\
&\delta_2(x_1[k_0])=w_2\beta_2(x_2[k_0]),\\
&w_2=|\rho_2|\tau(1+\lambda_{21})+1/\gamma_2,
\end{aligned}\right.
\end{equation}
and $\lambda_{21},...,\lambda_{25}$ are positive real numbers.

Define 
\begin{align}
x=\left(\begin{matrix}
x_1\\ x_2
\end{matrix}\right),\;
r&=\left(\begin{matrix}
r_1\\ r_2
\end{matrix}\right),\;
y=\left(\begin{matrix}
y_1\\ \tilde y_2
\end{matrix}\right),\;
u=\left(\begin{matrix}
u_1\\ \tilde u_2
\end{matrix}\right),\label{feed2}\\
V(x)&=\frac{1}{\tau}V_1(x_1)+\frac{1}{\tau}V_2(x_2).\label{feed3}
\end{align}
The passivity property of the system $T_\tau(P)\times_{\mathcal{F}}\tilde T_\tau(\S)$ in Fig. \ref{figfe} (a) is given by the following theorem whose proof can be obtained by using Lemma \ref{lemfeedind}.
\begin{theorem}\label{thmfeedback}
	Suppose that $P$ (resp. $\Sigma$) satisfies Assumption \ref{assadd} and Assumption \ref{assgain} with constant $\gamma_1$ and function $\beta_1$ (resp. with constant $\gamma_2$ and function $\beta_2$). Suppose also that $P$ (resp. $\Sigma$) is  IF-OFP$(\nu_1,\rho_1)$ (resp. IF-OFP$(\nu_2,\rho_2)$) with a continuous, positive definite, radially unbounded storage function  $V_1$ (resp. $V_2$). 
	Then the system $T_\tau(P)\times_{\mathcal{F}}\tilde T_\tau(\S)$ in Fig. \ref{figfe} (a) satisfies the following passivity  inequality for any $k_0,k_1\in\Z_{\geq 0},k_0<k_1$ and any admissible input $r[k]\in\R^{2m}$:
	\begin{align}
	&V(x[k_1])-V(x[k_0])\leq \delta(x[k_0])+\sum_{k=k_0}^{k_1-1} \big(r[k]^\top y[k]\big.\nonumber\\ 
	&\quad\quad\big.-\hat \nu |r[k]|^2-\hat \rho |y[k]|^2+\tilde \delta_2\big) \label{con:ineq3}
	\end{align}
	where $\delta(x[k_0])=\delta_1(x_1[k_0])+\delta_2(x_2[k_0])$, 
	\begin{equation}\label{composedind}
	\left\{
	\begin{aligned}
	&\hat\nu<\min\{\nu_1',\tilde\nu_2\},\\
	&\hat\rho\leq \min\{\rho_1'-\frac{\hat\nu\tilde\nu_2}{\tilde\nu_2-\hat{\nu}},\tilde\rho_2-\frac{\hat \nu\nu_1'}{\nu_1'-\hat{\nu}}\},
	\end{aligned}\right.
	\end{equation}
	and  $\delta_1,\delta_2$, $\tilde \delta_2$, $\nu_1',\rho_1'$, $\tilde \nu_2,\tilde \rho_2$ are given in \eqref{comind1} and \eqref{comind2}.
\end{theorem}



The following corollary shows the passivity inequality $T_\tau(P)\times_{\mathcal{F}}\tilde T_\tau(\S)$ satisfies for a special case when $r[k]\equiv{\bf 0}$.
\begin{corollary}\label{corfeed}
	Suppose that conditions in Theorem \ref{thmfeedback} hold.
	When $r[k]\equiv {\bf 0}$, the system $T_\tau(P)\times_{\mathcal{F}}\tilde T_\tau(\S)$ in Fig. \ref{figfe} (a) satisfies the following inequality 
	\begin{align}
	&V(x[k_1])\!-\!V(x[k_0])\leq \delta(x[k_0])\!+\!\sum_{k=k_0}^{k_1-1} \big[\!-\!(\tilde \nu_2+\rho_1')|y_1[k]|^2\big.\nonumber\\ 
	&\big.-( \nu_1'\!+\!\tilde\rho_2) |\tilde y_2[k]|^2+\tilde \delta_2\big],\;\forall k_0,k_1\in\Z_{\geq 0},k_0<k_1,\label{con:ineq5}
	\end{align}
	where $\delta(x[k_0])=\delta_1(x_1[k_0])+\delta_2(x_2[k_0])$, and  $\delta_1,\delta_2,\tilde \delta_2$, $\nu_1',\rho_1'$, $\tilde \nu_2,\tilde \rho_2$ are given in \eqref{comind1} and \eqref{comind2}.
\end{corollary}

{\color{black}The particular case when $r[k]\equiv{\bf 0}$ requires less restrictive conditions for the closed-loop system to have positive passivity indices. 
By Remark \ref{remark8}, when $r[k]\equiv{\bf 0}$, if $\nu_2'+\rho_1'>0$ and $\nu_1'+\rho_2'>0$,} then it is always possible to choose $\lambda_{22},...,\lambda_{25}$ large enough and $\mu_1,\mu_2$ small enough such that $\tilde \nu_2+\rho_1'>0$, $\nu_1'+\tilde\rho_2>0$ and $\tilde \delta_2$ arbitrarily small. {\color{black}On the contrary, \eqref{composedind} implies that $\hat \rho$ might be negative even if $\nu_1',\nu_2',\rho_1',\rho_2'$ are all positive. }

\subsection{Ultimate Boundedness Of the Closed-Loop System}

{\color{black}The following theorem shows conditions under which the state of $T_\tau(P)\times_{\mathcal{F}}\tilde T_\tau(\S)$ is ultimately bounded. The proof of the theorem is given in Appendix \ref{proofthm3}.
\begin{theorem}\label{thm3}
	Consider the system $T_\tau(P)\times_{\mathcal{F}}\tilde T_\tau(\S)$ in Fig. \ref{figfe} (a). Suppose that 1) $T_{\tau}(P)$ and $\tilde T_{\tau}(\S)$ satisfy passivity inequalities \eqref{con:ineq1} and \eqref{con:ineq2}, respectively; 2) $T_\tau(P)$ is $N_1$-step SD that satisfies \eqref{det1} and $T_\tau(\S)$ is $N_2$-step SD that satisfies  \eqref{det2}; 3) $\hat \rho$ shown in \eqref{composedind} can be chosen as $\hat \rho>0$; 4) there exist a compact set $\mathcal{X}\subset \R^{n_1+n_2}$ and a radially unbounded, positive definite function $\kappa: \R^{n_1+n_2}\rightarrow\R$ such that $\mathcal{X}\supseteq \mathcal{D}_1$  and $\eta_2p(x)-\delta(x)\geq \kappa(x)$ for all $x\in\mathcal{X}$ where $\mathcal{D}_1:=\{z\mid V(z)\leq \max\{V(x[0]),V(x[1]),...,V(x[N]),d_1+d_2+d_4\}\}$ is a compact set,  $V$ is defined in \eqref{feed3}, $d_1:=(N+1)[(\eta_1+\vartheta\eta_2) \|r\|^2+\tilde\delta_2]$,
	$d_2:=m(1-\vartheta)(N_2+1)(\vartheta_2\mu_1^2+\mu_2^2)$,
	$d_3$ is a positive real number,
	$d_4:=\max_{z\in\mathcal{C}}V(z)$, 
	$\mathcal{C}:=\{z\mid p(z)\leq d_1+d_2+ d_3\}$,  $p(x):=(1-\vartheta)(p_1(x_1)+\frac{1}{2}p_2(x_2))$, $\eta_1:=\frac{1}{4\l}-\hat\nu> 0$, $\eta_2:=\hat\rho-\l>0$, $\l$ is a number satisfying $0<\l<\hat\rho$, $\vartheta$ is defined in \eqref{vartheta}, $\delta$ and $\tilde\delta_2$ are given in Theorem \ref{thmfeedback}, and $N:=\max\{N_1,N_2\}$.
	Then,\\
	i) $x[k]\in\mathcal{D}_1$ for any $k\in\Z_{\geq 0}$;\\ 
	ii) there exists $K\in\Z_{\geq 0}$ such that $x[k]\in\mathcal{D}_2$ for all $k\geq K$ where  $\mathcal{D}_2:=\{x\mid V(x)\leq d_1+d_2+d_4\}$  is a compact set.
\end{theorem}}

When  $r[k]\equiv{\bf 0}$, conditions in Theorem \ref{thm3} can be simplified, as shown by the following corollary. 
{\color{black}\begin{corollary}\label{cor4}
	Consider the system $T_\tau(P)\times_{\mathcal{F}}\tilde T_\tau(\S)$ in Fig. \ref{figfe} (a) where $r[k]\equiv{\bf 0}$. Suppose that conditions 1)-2) in Theorem \ref{thm3} hold, and 3) $\tilde \nu_2+\rho_1'>0$,  $\nu_1'+\tilde\rho_2>0$; 4) there exist a compact set $\mathcal{X}\subset \R^{n_1+n_2}$ and a radially unbounded, positive definite function $\kappa: \R^{n_1+n_2}\rightarrow\R$ such that $\mathcal{X}\supseteq \mathcal{D}_1$  and $p(x)-\delta(x)\geq \kappa(x)$ for all $x\in\mathcal{X}$ where $\mathcal{D}_1:=\{z\mid V(z)\leq \max\{V(x[0]),V(x[1]),...,V(x[N]),d_1+d_2+d_4\}\}$,  $V$ is defined in \eqref{feed3}, $d_1:=(N+1)\tilde\delta_2$,
	$d_2:=m(N_2+1)(\vartheta_2\mu_1^2+\mu_2^2)$,
	$d_3$ is a positive real number,
	$d_4:=\max_{z\in\mathcal{C}}V(z)$, 
	$\mathcal{C}:=\{z\mid p(z)\leq d_1+d_2+ d_3\}$ and $p(x):=p_1(x_1)+\frac{1}{2}p_2(x_2)$, $\delta$ and $\tilde\delta_2$ are given in Theorem \ref{thmfeedback}, and $N:=\max\{N_1,N_2\}$. Then, i) $x[k]\in\mathcal{D}_1$ for any $k\in\Z_{\geq 0}$; ii)  there exists $K\in\Z_{\geq 0}$ such that $x[k]\in\mathcal{D}_2$ for all $k\geq K$ where  $\mathcal{D}_2:=\{x\mid V(x)\leq d_1+d_2+d_4\}$  is a compact set.
\end{corollary}}

The proof of Corollary \ref{cor4} is given in Appendix \ref{proofthm4}.
Corollary \ref{cor4} shows that the state of  $T_\tau(P)\times_{\mathcal{F}}\tilde T_\tau(\S)$ with $r[k]\equiv{\bf 0}$ is ultimately bounded, if the shortage of passivity of one component can be compensated for by the excess of passivity of another component in the feedback connection, and in addition,  each component is SD. 

{\color{black}When the dynamical models of a plant and a controller are known exactly, it might be easier to analyze the closed-loop system using these models directly. However, it is almost impossible to obtain an exact model of a plant or a controller. The passivity-based analysis shown above enables us to compute the passivity degradation from the original indices and derive the property of the closed-loop system (under sampling and quantization) without knowing the exact model {\color{black}(e.g., see \cite{wu2013experimentally,romer2018some,antsaklis2013control})}. Therefore, the passivity-based analysis is robust because it holds for a family of plants and controllers. }
Although the ultimate bounds given in Theorem \ref{thm3} and Corollary \ref{cor4} may not be  tight, we point out that the ultimate bound given in Theorem \ref{thm3}  will decrease as $\|r\|,\mu_1,\mu_2$ decrease, and the ultimate bound given in Corollary \ref{cor4} can be made arbitrarily small when $\mu_1,\mu_2$ are chosen small enough.

\begin{remark}
	In \cite{xupasscdc17}, modifying the passivity indices of a system by simply adding a feedforward loop and/or a feedback loop 
	was discussed, where the achievable bounds of the modified passivity indices were given explicitly.  When condition 4) in Theorem \ref{thm3} (or Corollary \ref{cor4}) is not satisfied, the results in \cite{xupasscdc17} may be used to modify the passivity indices of the feedback-connected system. 
\end{remark}

\begin{remark}
	The emulation-based design is a general framework for the controller design under time sampling: a  controller is designed in the continuous time domain at first, and then sampled and implemented using a sampler and hold device so that certain property can be preserved under the time sampling (see \cite{LailaSampled06,LailaECC02,nesic2001sampled,nesic2009explicit} and references therein). The notion of $(V,w)$-dissipativity was used to characterize the property of the system for preservation.
	There are two main differences between the emulation-based approach and the method used here: 1) the sampling time in the emulation-based approach is an implicit parameter that needs to be chosen sufficiently small but whose explicit value is hard to compute in general  (see \cite{nesic2009explicit} for more details), while the sampling time in the results above is an explicit parameter whose effect on the degradation of the passivity indices is explicitly known; 2) the input and output quantizations are not considered in the emulation-based approach. In the next section, we will show that the symbolic control implementation can be also studied in our framework.
\end{remark}
\section{Analysis of the Approximate Bisimulation-Based Control Implementation}\label{sec:symbolic}

\begin{figure}[!h]
	\centering
	\includegraphics[width=\linewidth]{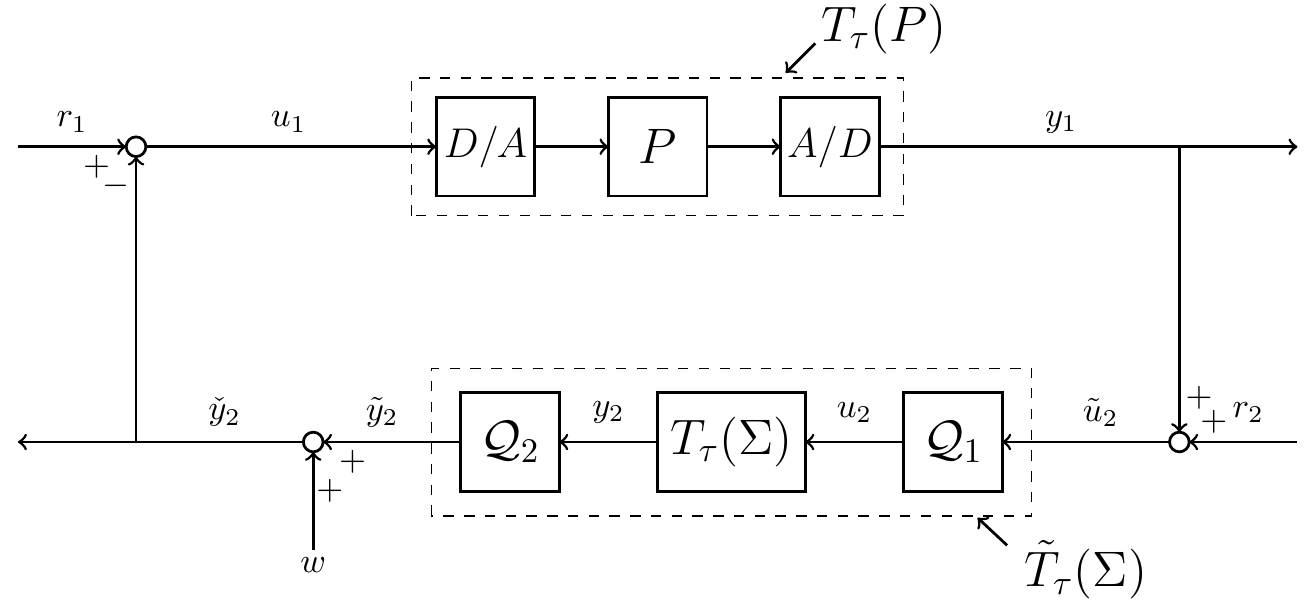}
	\caption{ {\color{black} Setup of the closed-loop system $T_\tau(P)\times_{\mathcal{F}}^w\tilde T_{\tau}(\S)$, which is $T_\tau(P)\times_{\mathcal{F}}\tilde T_{\tau}(\S)$ with a bounded disturbance $w$ adding to $\tilde y_2$.} }\label{figfe3}
\end{figure}

In this section, we study state boundedness of the system $T_\tau(P)\times_{\mathcal{F}}\tilde T_{\tau\mu\eta}(\S)$ shown in Fig. \ref{figfe} (b) where a symbolic controller $T_{\tau\mu\eta}(\S)$ approximately bisimilar to $\S$ is implemented. 
To that end, we introduce an auxiliary configuration, denoted by $T_\tau(P)\times_{\mathcal{F}}^w\tilde T_{\tau}(\S)$, as shown in Fig. \ref{figfe3} where an external bounded disturbance $w$ is added to $\tilde y_2$ in the system $T_\tau(P)\times_{\mathcal{F}}\tilde T_{\tau}(\S)$. In what follows, we will first show the passivity property and state boundedness of $T_\tau(P)\times_{\mathcal{F}}^w\tilde T_{\tau}(\S)$ by assuming that $w$ is bounded, and based on that, we will show the state boundedness of $T_\tau(P)\times_{\mathcal{F}}\tilde T_{\tau\mu\eta}(\S)$ by showing that there is a particular choice of $w$ that relates the auxiliary configuration and the setup $T_\tau(P)\times_{\mathcal{F}}\tilde T_{\tau\mu\eta}(\S)$ of interest. 

Denote the states of $T_{\tau}(P)$,  $T_{\tau}(\S)$ and $T_{\tau\mu\eta}(\S)$ by $x_1[k]$, $x_2[k]$ and $x^s_2[k]$, respectively (``s'' stands for ``symbolic'').  
Suppose that $P$ (resp. $\Sigma$) satisfies Assumption \ref{assadd} with $y_1(x_1,u_1)=h_{11}(x_1)+h_{12}(u_1)$ (resp. with $y_2(x_2,u_2)=h_{21}(x_2)+h_{22}(u_2)$) and Assumption \ref{assgain} with constant $\gamma_1$ and function $\beta_1$ (resp. constant $\gamma_2$ and function $\beta_2$). Suppose also that $P$ (resp. $\Sigma$) is  IF-OFP$(\nu_1,\rho_1)$ (resp. IF-OFP$(\nu_2,\rho_2)$) with a continuous, positive definite, radially unbounded storage function  $V_1$ (resp. $V_2$). Furthermore, suppose that there exist positive numbers $\epsilon,\mu,\eta>0$ such that $T_\tau(\S)\cong^{(\e,\mu)} T_{\tau\mu\eta}(\Sigma)$ holds where $\mu=\mu_1$ {\color{black}(recall that $\mu_1$ is the quantization level of the quantizer $\mathcal{Q}_1$)}, and there exists $L>0$ such that for all $z_1,z_2\in\R^n$,
\begin{align}\label{ineqL}
|h_{21}(z_1)-h_{21}(z_2)|\leq L|z_1-z_2|_\infty.
\end{align}
Define $x=\left(\begin{matrix}
x_1\\ x_2
\end{matrix}\right),\;
r=\left(\begin{matrix}
r_1\\ r_2
\end{matrix}\right),\;
y=\left(\begin{matrix}
y_1\\ \check y_2
\end{matrix}\right),\;
u=\left(\begin{matrix}
u_1\\ \tilde u_2
\end{matrix}\right)$ 
for signals depicted in Fig. \ref{figfe3}.
Next, we show the passivity property of $T_\tau(P)\times_{\mathcal{F}}^w\tilde T_{\tau}(\S)$ by assuming that
\begin{align}
|w[k]|\leq L\epsilon + 2\sqrt{m}\mu_2, \;\forall k\in\Z_{\geq 0},\label{wbound}
\end{align}
where $\mu_2$ is the quantization level of quantizer $\mathcal{Q}_2$.
Recall that for any $k \in\Z_{\geq 0}$, $|u_2[k]| \leq  |\tilde u_2[k]|$, $|u_2[k]-\tilde u_2[k]|\leq \sqrt{m}\mu_1$ and $|\tilde y_2[k]-y_2[k]|\leq \sqrt{m}\mu_2$. In addition, since $\check y_2[k] = w[k]+\tilde y_2[k]$, 
it follows that $|\check y_2[k]-y_2[k]|\leq |\check y_2[k]-\tilde y_2[k]|+|\tilde y_2[k]-y_2[k]|\leq L\epsilon+3\sqrt{m}\mu_2,\forall k \in\Z_{\geq 0}.$
Similar to Theorem \ref{thm2}, it is easy to show that the following inequality holds  for any $k_0,k_1\in\Z_{\geq 0},k_0<k_1$, and any admissible inputs $\tilde u_2\in U_1$:
\begin{align}
&\hat V_2(x_2[k_1])-\hat V_2(x_2[k_0]) \leq \delta_2+\sum_{k=k_0}^{k_1-1} \Big(\tilde u_2[k]^\top \check y_2[k]\big.\nonumber\\
&\quad\quad\big.-\tilde\nu_2|\tilde u_2[k]|^2-\tilde\rho_2|\check y_2[k]|^2+\tilde \delta_2\Big)\label{ineqsym}
\end{align}
where $\hat V_2(x_2)=\frac{1}{\tau}V_2(x_2)$, $\delta_2$, $\tilde \nu_2$, $\tilde \rho_2$ are given by the equations in \eqref{comind2}, respectively, and
\begin{align}\label{indexsym}
\tilde \delta_2=&[|\rho_2'|(1+\l_{22})+\l_{24}](L\epsilon + 3\sqrt{m}\mu_2)^2\nonumber\\
&\quad+[|\nu'_2|(1+\l_{23})+\l_{25}]m\mu_1^2,
\end{align}
with $\nu'_2$, $\rho'_2$ given in \eqref{comind2}. 

As discussed in Subsection \ref{subpassineq},  $\tilde T_{\tau}(P)$ satisfies \eqref{con:ineq1} where $\nu_1',\rho_1',\delta_1$ are given in \eqref{comind1}. Similar to Theorem \ref{thmfeedback}, it is easy to show that  $T_\tau(P)\times_{\mathcal{F}}^w\tilde T_{\tau}(\S)$ shown in Fig. \ref{figfe3} satisfies inequality \eqref{con:ineq3} for $k_0,k_1\in\Z_{\geq 0},k_0<k_1$ and any admissible input $r[k]\in\R^{2m}$
where  $V(x)=\frac{1}{\tau}V_1(x_1)+\frac{1}{\tau}V_2(x_2)$, $\hat\nu,\hat\rho$ satisfy \eqref{composedind}, 
$\delta(x[0])=\delta_1(x_1[0])+\delta_2(x_2[0])$, and  $\delta_1,\delta_2,\tilde \delta_2,\nu_1',\rho_1',\tilde \nu_2,\tilde \rho_2$ are those in \eqref{comind1} and \eqref{ineqsym}.

Now let us consider the SD property of $T_\tau(P)\times_{\mathcal{F}}^w\tilde T_{\tau}(\S)$ similar to the discussion in preceding sections. Suppose that $T_\tau(P)$ and $T_\tau(\S)$ are both SD such that inequalities \eqref{det1} and \eqref{det2} hold. 
Note that $|\tilde u_2|^2\geq \frac{1}{2}|u_2|^2-|\tilde u_2-u_2|^2\geq\frac{1}{2}|u_2|^2-m\mu_1^2$ and $|\check y_2|^2\geq \frac{1}{2}|y_2|^2-|\check y_2-y_2|^2\geq\frac{1}{2}|y_2|^2-(L\epsilon+3\sqrt{m}\mu_2)^2$.
Then, we have $\vartheta_2|\tilde u_2[k]|^2+|\check y_2[k]|^2
\geq \frac{1}{2}(\vartheta_2|u_2[k]|^2+|y_2[k]|^2)-\vartheta_2m\mu_1^2-(L\epsilon+3\sqrt{m}\mu_2)^2$, 
which implies that $\sum_{k=k_0}^{k_0+N_2}\vartheta_2|\tilde u_2[k]|^2+|\check y_2[k]|^2\geq \frac{1}{2}p_2(x_2[k_0]) -(N_2+1)[\vartheta_2m\mu_1^2+(L\epsilon+3\sqrt{m}\mu_2)^2].$
Then, similar to the proof of Lemma \ref{lemstrongdet}, it is easy to show that for any $x[k_0]$ and any $r[k]$, we have 
\begin{align*}
&\sum_{k=k_0}^{k_0+N} \vartheta |r[k]|^2+|y[k]|^2\geq p(x[k_0])\\
&\quad\quad\quad-(N_2+1)[\vartheta_2m\mu_1^2+(L\epsilon+3\sqrt{m}\mu_2)^2]
\end{align*}
where $N,p,\vartheta$ are those given in Theorem \ref{thm3}.

{\color{black}With notations above, we have the following Lemma \ref{thm5} for the system $T_\tau(P)\times_{\mathcal{F}}^w\tilde T_{\tau}(\S)$ shown in Fig. \ref{figfe3}. The proof of Lemma \ref{thm5} is similar to that of Theorem \ref{thm3} and is omitted due to the space limitation.}
{\color{black}\begin{lemma}\label{thm5}
	Consider the system $T_\tau(P)\times_{\mathcal{F}}^w\tilde T_{\tau}(\S)$ shown in Fig. \ref{figfe3} where $w[k]$ is bounded and satisfies \eqref{wbound} with $L$ given in \eqref{ineqL} and $\epsilon,\mu_2>0$. Suppose that inequalities \eqref{con:ineq1} and \eqref{ineqsym} hold and conditions 2)-4) in Theorem \ref{thm3} hold.
	Then, there exist compact sets $\mathcal{D}_1,\mathcal{D}_2$ and a constant $K\in\Z_{\geq 0}$ such that for any $w[k]$ that satisfies \eqref{wbound}, $(x_1[k],x_2^s[k])\in\mathcal{D}_1$ for  $k\in\Z_{\geq 0}$ and $(x_1[k],x_2^s[k])\in\mathcal{D}_2$ for  $k\geq K$. 
\end{lemma}


Using Lemma \ref{thm5}, we have the following results that provide conditions for the state boundedness of  $T_\tau(P)\times_{\mathcal{F}}\tilde T_{\tau\mu\eta}(\S)$ shown in Fig. \ref{figfe} (b). The proof of Theorem \ref{corsym} is given in Appendix \ref{proofcor6}.
\begin{theorem}\label{corsym}
	Consider the system $T_\tau(P)\times_{\mathcal{F}}\tilde T_{\tau\mu\eta}(\S)$ shown in Fig. \ref{figfe} (b). Suppose that all the conditions of Lemma \ref{thm5} hold, and the initial conditions  $x_2[0]$ and $x_2^s[0]$ satisfy  $|x_2[0]-x_2^s[0]|_\infty\leq \epsilon$.
	Then,  there exist compact sets $\mathcal{D}_1,\mathcal{D}_2$ and a constant $K\in\Z_{\geq 0}$ such that $(x_1[k],x_2^s[k])\in\mathcal{D}_1$ for  $k\in\Z_{\geq 0}$ and $(x_1[k],x_2^s[k])\in\mathcal{D}_2$ for  $k\geq K$.
\end{theorem}

\begin{corollary}\label{cor6}
	Consider the system $T_\tau(P)\times_{\mathcal{F}}\tilde T_{\tau\mu\eta}(\S)$ shown in Fig. \ref{figfe} (b) where 
	 $r[k]\equiv {\bf 0}$. Suppose that 1) inequalities \eqref{con:ineq1} and \eqref{ineqsym} hold, 2) $|x_2[0]-x_2^s[0]|_\infty\leq \epsilon$, 3) condition 2) of Theorem \ref{thm3} holds, and 4) conditions 3)-4) in Corollary \ref{cor4} hold. Then,  there exist compact sets $\mathcal{D}_1,\mathcal{D}_2$ and a constant $K\in\Z_{\geq 0}$ such that $(x_1[k],x_2^s[k])\in\mathcal{D}_1$ for  $k\in\Z_{\geq 0}$ and $(x_1[k],x_2^s[k])\in\mathcal{D}_2$ for  $k\geq K$.
\end{corollary}}
Similar to the discussion in Section \ref{sec:connect}, the bounds of $x_1,x_2^s$  can be made arbitrarily small by letting the parameters $\mu_1,\mu_2,\eta,\epsilon$ small enough. Furthermore, the bounds of $x_2^s$ can be used to determine how the state space and input space of $T_{\tau\m\eta}(\S)$ can be chosen as compact sets; the resulting $T_{\tau\m\eta}(\S)$ with bounded state space and input space is a finite transition system that can be implemented with finite precision. 

Above result formalizes the intuition that when a controller guaranteeing stability in some robust way is replaced by its symbolic bisimilar version, the feedback-connected system should be ``somewhat stable". The key challenge to keep in mind when replacing the controller in the feedback loop is the fact that the internal signals driving the plant changes therefore, within the feedback-loop, symbolic model and the actual controller can be driven by totally different inputs even if they are initialized with the same initial conditions. The external bounded signal $w$ we introduce in our analysis captures this additional robustness required. Moreover, passivity indices provide us a way to explicitly compute global and ultimate bounds on the states of the closed-loop system with the symbolic controller implementation. 

\begin{remark}
	There is a trade-off between the passivity degradation and the construction of the approximately bisimilar model $T_{\tau\m\eta}(\S)$: on one hand, a smaller sampling time $\tau$ can result in a smaller passivity degradation; on the other hand, given a precision $\epsilon$, the sampling time
	$\tau$ should be chosen large enough in order to make \eqref{eqnbisi} holds and  $T_\tau(\S)\cong ^{(\e,\mu)}T_{\tau\m\eta}(\S)$. 
\end{remark}

{\color{black}	\begin{example}\label{ex5}
	Consider the following nonlinear system that is adopted from Example 1 of \cite{xia2017approx}: 
	\begin{align*}
	P: \begin{cases}
	\dot{x}_1=-0.7x_1-0.2x_1^3-0.5x_2+0.4u_1\\
	\dot{x}_2=0.5x_1-0.3x_2^3+0.5u_2\\
	y_1=0.4x_1\\
	y_2=0.5x_2
	\end{cases}
	\end{align*}
	System $P$ is passive with a storage function $V(x)=\frac{1}{2}(x_1^2+x_2^2)$ since  $\dot{V}\leq u^\top y$, which implies $P$ is IF-OFP$(0,0)$. Assumption 1 holds for $P$ since $h_1(x)=(0.4x_1,0.5x_2)^\top$, and Assumption 2 holds for $h_1$ with $\gamma=0.8189$ since \eqref{gainlya} can be verified using sum-of-squares optimization. 
	Consider the configuration of Fig. \ref{figfe} (b) where $P$ is the plant above and $\Sigma$ is the controller  in Example \ref{ex1}.
	Suppose that  $r[k]={\bf 0}$  and choose $\mu_1=\mu_2=0.01$. 
	Given a precision $\epsilon=0.25$, parameters $\tau=0.3$, $\mu=\mu_1$, $\eta=0.1$ can be chosen such that $T_\tau(\S)\cong ^{(\e,\mu)}T_{\tau\m\eta}(\S)$. We can verify that $T_\tau(P)$ is $0$-step SD with $\vartheta_1=0$ and $p_1(x)=0.16x_1^2+0.25x_2^2$, and $T_\tau(\S)$ is $0$-step SD with $\vartheta_2=0$ and $p_2(x)=0.13x_1^2-0.03x_1x_2+0.1125x_2^2$. We can also verify that conditions of Corollary \ref{cor6} hold, and thus $x_1[k],x_2^s[k]$ are ultimately bounded. 
	To demonstrate the effect of state quantization of $T_{\tau\m\eta}(\S)$, we also consider using quantization precisions $\eta=0.05$ and $\eta=0.01$, both of which also guarantee  $T_\tau(\S)\cong ^{(\e,\mu)}T_{\tau\m\eta}(\S)$.
	In Fig. \ref{figureex5}, the trajectories of $x_1[k]$, $x_2^s[k]$ are shown for three different choices of $\eta$ where the initial states are all $x_1[0]=(-0.7,-2)^\top$, $x_2[0]=(1.5,-1.6)^\top$. 
	From these figures, it can be observed  that the state trajectories are eventually bounded, and the smaller $\eta$ (i.e., the state quantization precision) is, the smaller the ultimate bound would be.  
	\begin{figure}[!h]
		\centering
		\begin{subfigure}[t]{0.5\textwidth}
			\centering\caption{$\eta=0.1$}
			\includegraphics[width=0.46\textwidth]{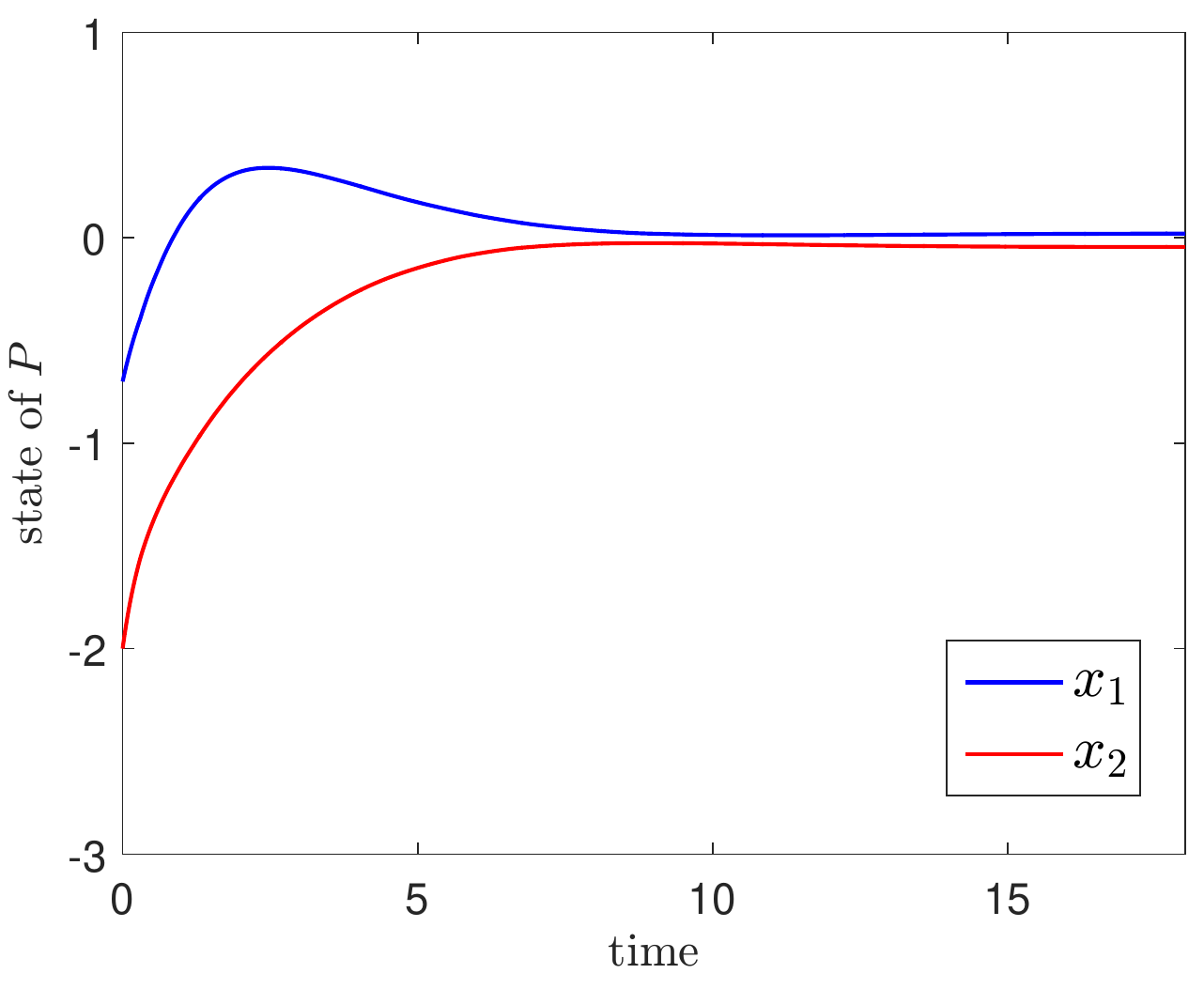}
			\includegraphics[width=0.46\textwidth]{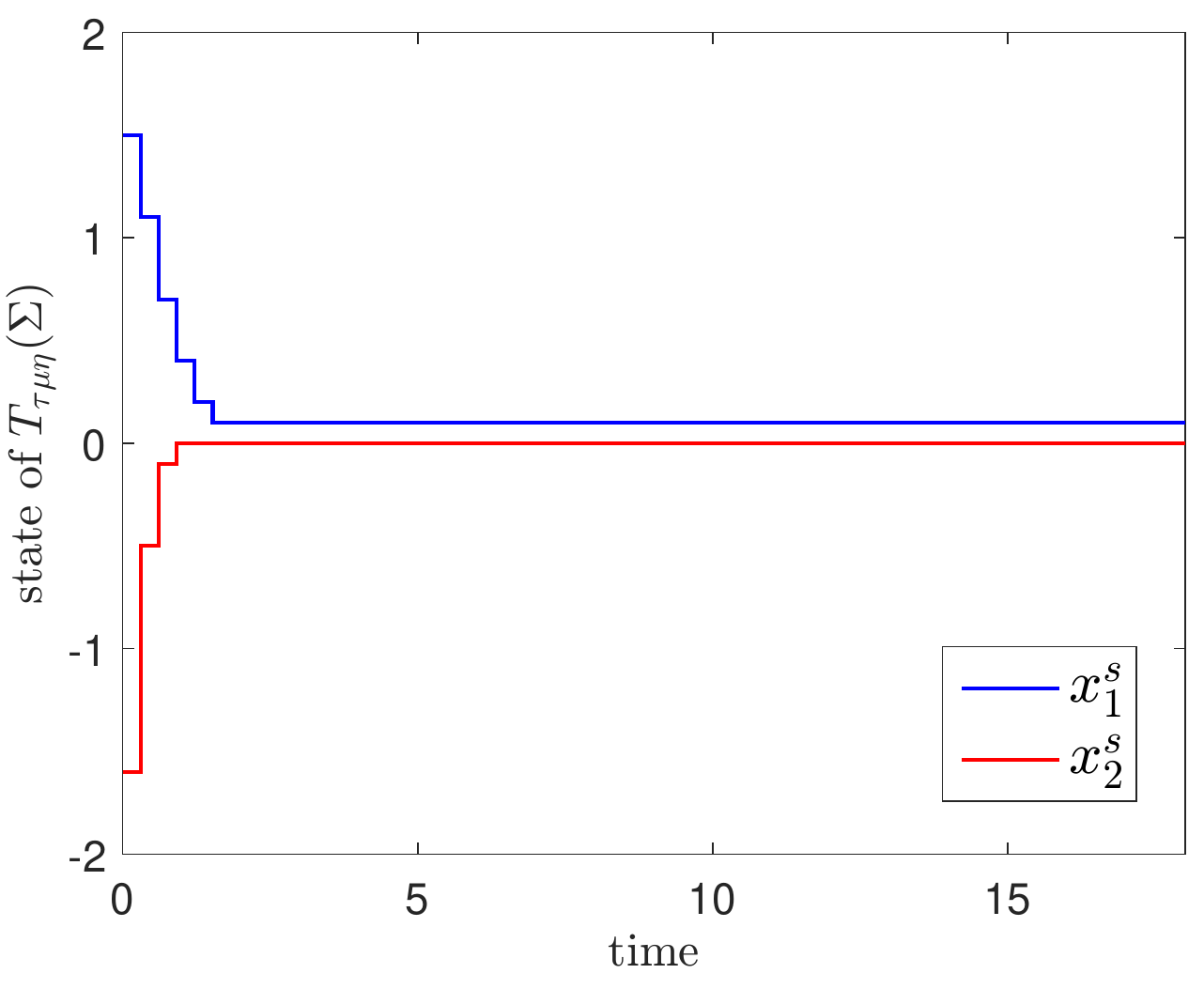}					
		\end{subfigure}
		\begin{subfigure}[t]{0.5\textwidth}
			\centering\caption{$\eta=0.05$}
			\includegraphics[width=0.46\textwidth]{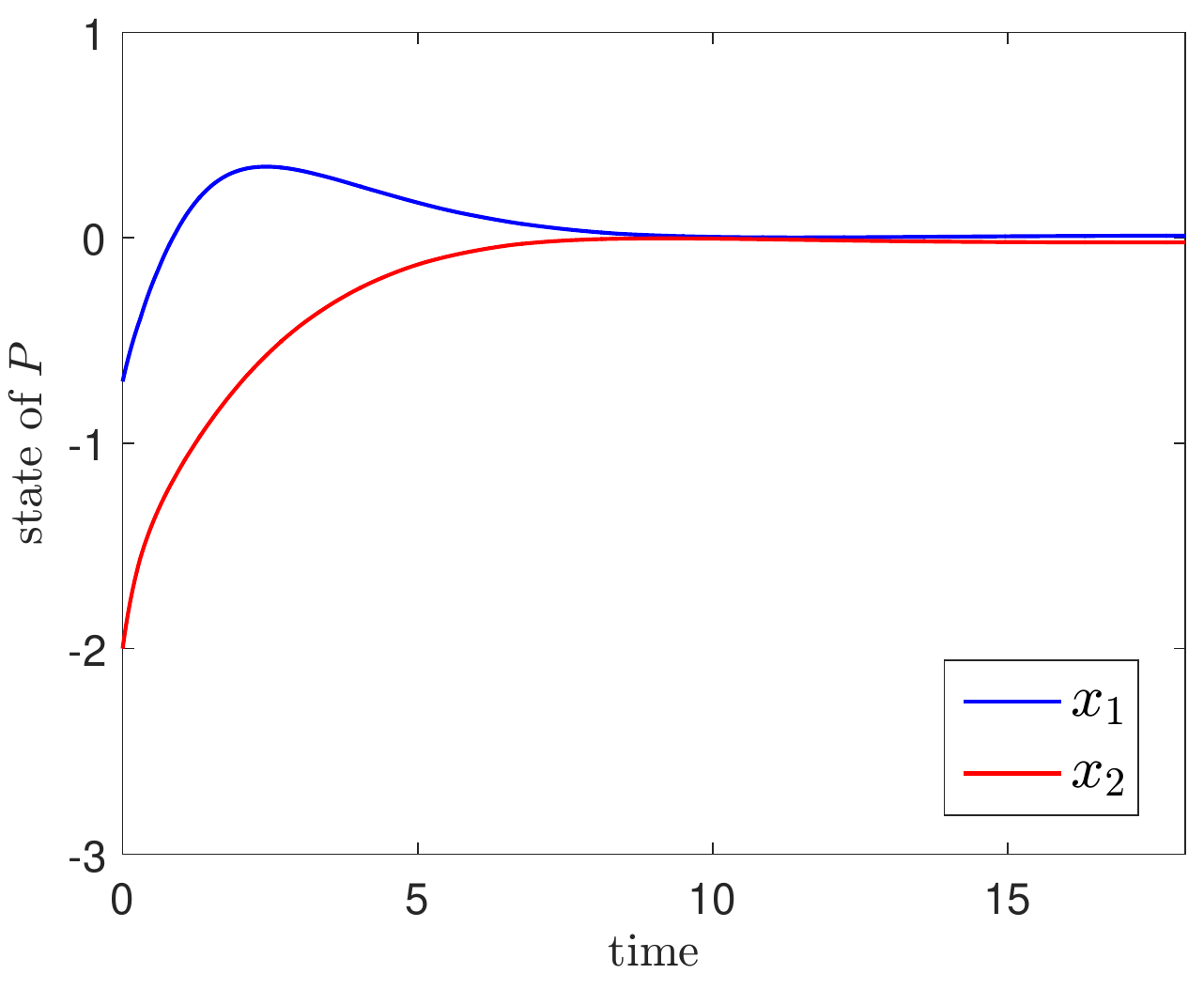}
			\includegraphics[width=0.46\textwidth]{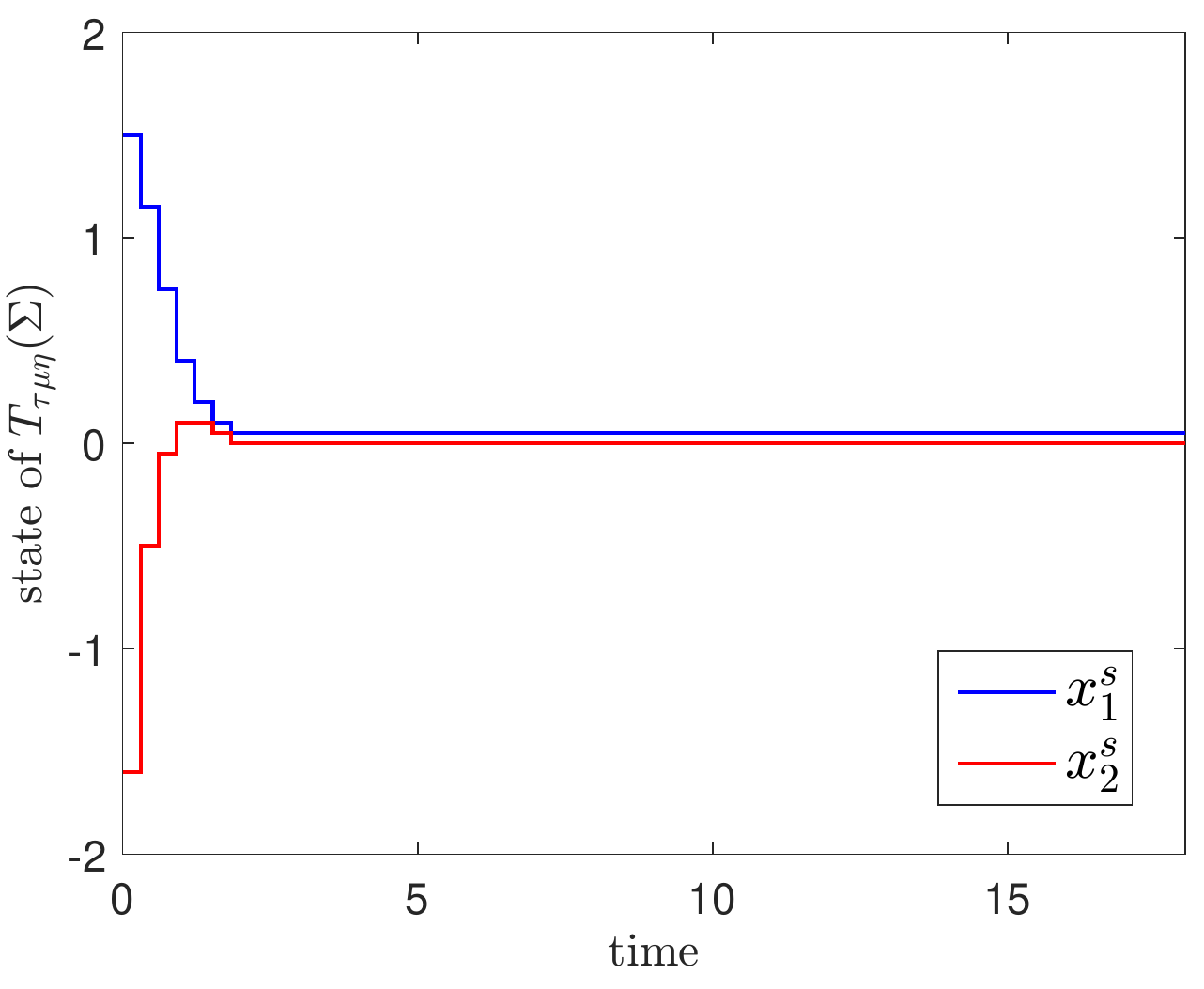}	
		\end{subfigure}
		\begin{subfigure}[t]{0.5\textwidth}
			\centering\caption{$\eta=0.01$}
			\includegraphics[width=0.46\textwidth]{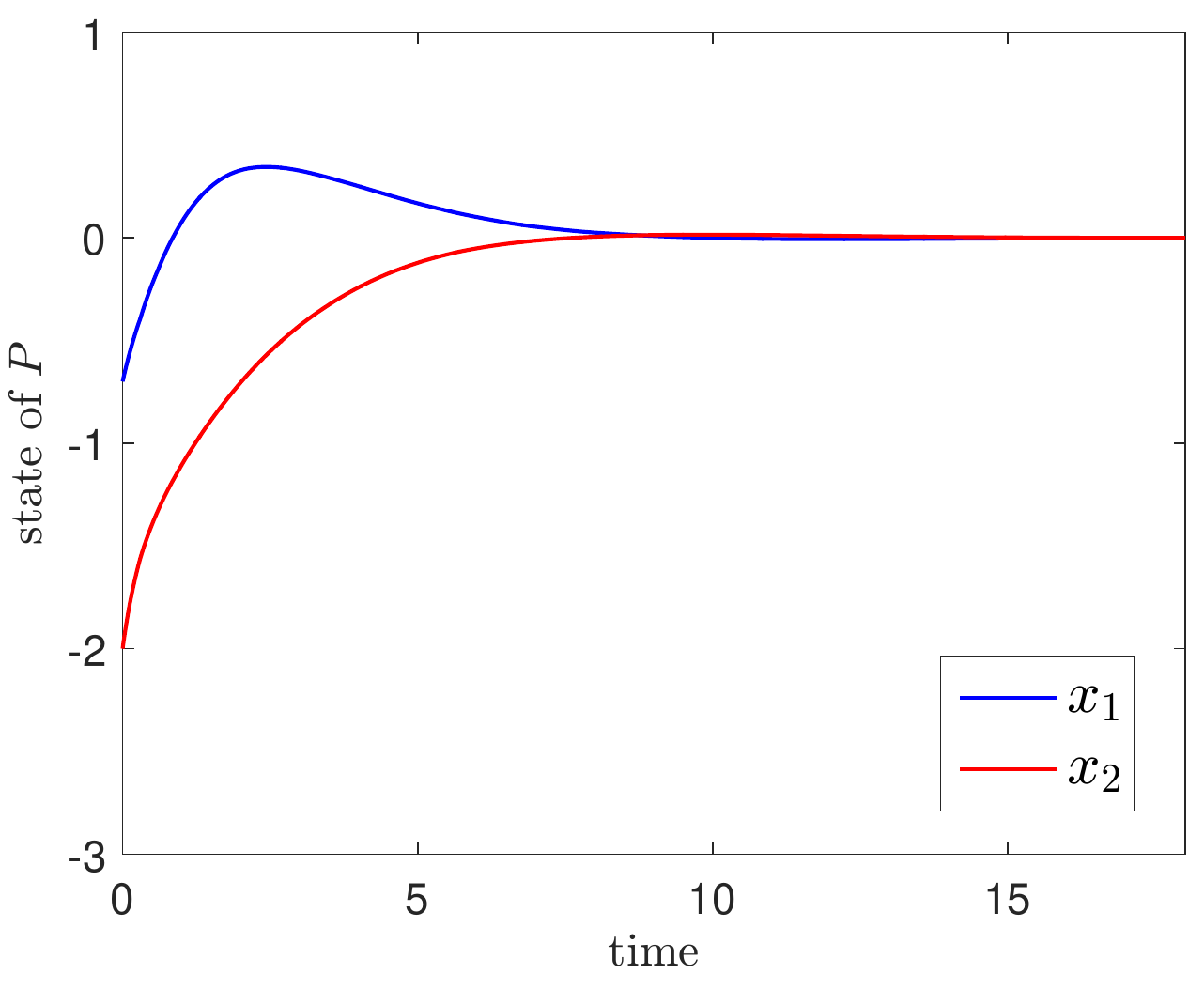}
			\includegraphics[width=0.46\textwidth]{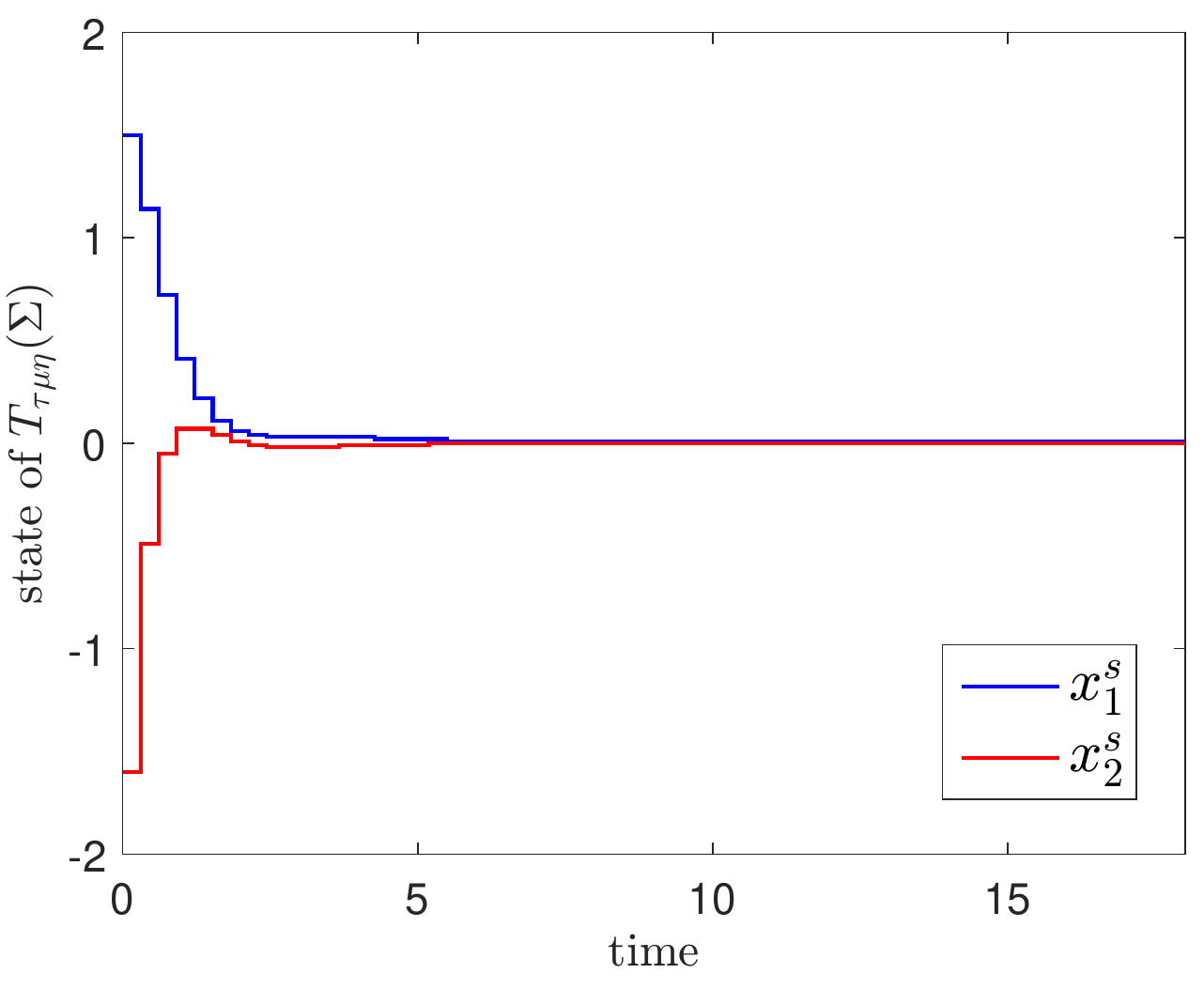}		
		\end{subfigure}
		\caption{State trajectories of $P$ and $T_{\tau\m\eta}(\S)$ using three $\eta$. }\label{figureex5}
	\end{figure}
\end{example}}




\section{Conclusions}\label{sec:conclusions}
In this paper, we analyzed the passivity degradation from a continuous-time control system to its sampled, input and output quantized model, and used these results to analyze the feedback-connected system when a continuous controller is implemented via sampling and quantization. 
We also considered bisimulation-based symbolic implementations, where in addition to inputs and outputs, the internal states of the controller is also quantized, therefore the controller has a fully discrete representation.
We derived conditions under which the closed-loop system with those control implementations is (ultimately) bounded, where explicit bounds are presented and ways to render these bounds arbitrarily small are discussed. Our method provides a novel perspective for the analysis of heterogeneous, compositional systems, and can be used in the analysis of potentially complex control systems with various implementation artifacts due to resource constraints. 
\bibliographystyle{abbrv}

\bibliography{./pass}

\section*{Appendix}
\begin{spacing}{0.2}
\appendix

\section{Proof of Lemma \ref{lemstrongdet}}\label{sub:strongdetec}
Clearly $0\leq \max\{2\vartheta_1(1-\vartheta),2\vartheta_2(1-\vartheta)\}\leq \vartheta <1$ and $p$ is a radially unbounded, positive definite function. 
Recall that $r_1[k]=y_2[k]+u_1[k]$, $r_2[k]=u_2[k]-y_1[k]$. For any $x[k_0]$ and any $r[k]$, 
\begin{align*}
&\sum_{k=k_0}^{k_0+N} \vartheta |r[k]|^2+|y[k]|^2\\
= &  \sum_{k=k_0}^{k_0+N}  (1-\vartheta)(|y_1[k]|^2+|y_2[k]|^2)+\frac{\vartheta}{2}(|u_1[k]|^2+|u_2[k]|^2)\\
&\quad\quad + \vartheta [(\frac{|u_1[k]|^2}{2}+2|y_2[k]|^2+2u_1[k]^\top y_2[k])]\\
&\quad\quad + \vartheta [(\frac{|u_2[k]|^2}{2}+2|y_1[k]|^2-2u_2[k]^\top y_1[k])]\\
\geq &  \sum_{k=k_0}^{k_0+N}  (1-\vartheta)(|y_1[k]|^2+|y_2[k]|^2)+\frac{\vartheta}{2}(|u_1[k]|^2+|u_2[k]|^2)\\
\geq & (1-\vartheta) \sum_{k=k_0}^{k_0+N}  |y_1[k]|^2+|y_2[k]|^2+\vartheta_1|u_1[k]|^2+\vartheta_2|u_2[k]|^2\\
\geq & p(x[k_0]).
\end{align*}
This completes the proof.
\hfill$\Box$
%
\section{Proof of Lemma \ref{lembound}}\label{sec:prf0}
Because $\nu\rho\leq \frac{1}{4}$ (see \cite{antsaklis2013control}), it holds that 
$\eta_1> \frac{1}{4\rho}-\nu\geq0$.
Because $u^\top y\leq \frac{1}{4\l}|u|^2+\l|y|^2$,
from \eqref{ieq:VSQP} we have
\begin{align}
& V(x[k+1])-V(x[k])\nonumber\\
\leq& \frac{1}{4\l}|u[k]|^2+\l|y[k]|^2-\nu |u[k]|^2-\rho |y[k]|^2+\delta\nonumber\\
=& \eta_1|u[k]|^2-\eta_2|y[k]|^2+\delta
\leq \eta_1\|u\|^2+\delta.\label{quasidetec}
\end{align}
Therefore, for any $k_0\in\Z_{\geq 0}$,
\begin{align}
& V(x[k_0+N+1])-V(x[k_0])\nonumber\\
\leq&  \sum_{k=k_0}^{k_0+N} \left(\eta_1|u[k]|^2-\eta_2|y[k]|^2+\delta\right)\nonumber\\
\leq& c_2 -\eta_2p(x[k_0]),\label{quasidetec2}
\end{align}
where the second inequality is from \eqref{strongdetec}. 

i). Since $V(x)$ is radially unbounded, $\mathcal{C}_1,\mathcal{D}_1$ are compact sets and $c_1<\infty$.
Clearly, $x[0]\in\mathcal{D}_1$.
For $k\in\{1,...,N,N+1\}$, \eqref{quasidetec} implies 
\begin{align}
V(x[k])&\leq V(x[0])+k(\eta_1\|u\|^2+\delta)\leq c_1+c_3.\label{quasidetec3}
\end{align}
Therefore, $x[k]\in\mathcal{D}_1$  for any $k\in\{1,...,N+1\}$.
Consider now $k\in\{N+2,...,2N+2\}$. If there exists some $k^*\in\{N+2,...,2N+2\}$ such that $x[k^*]\notin\mathcal{D}_1$, or equivalently, 
\begin{equation}\label{bound}
V(x[k^*])>c_1+c_3,
\end{equation}
then by (\ref{quasidetec}) we have
$V(x[k^*-N-1])\geq V(x[k^*])-(N+1)(\eta_1\|u\|^2+\delta)>c_1,$
implying that $p(x[k^*-N-1])> c_2/\eta_2$. Then,
\begin{align}
V(x[k^*])&\leq  V(x[k^*-N-1])+c_2 -\eta_2p(x[k^*-N-1])\nonumber\\
&\leq  V(x[k^*-N-1])\leq c_1+c_3\label{quasidetec4}
\end{align}
where the first inequality is from  \eqref{quasidetec2}, 
and the third inequality is from \eqref{quasidetec3}. Noting that \eqref{quasidetec4} contradicts with (\ref{bound}), we conclude that $x[k]\in\mathcal{D}_1$  for $k\in\{N+2,...,2N+2\}$. By induction, $x[k]\in\mathcal{D}_1$ for any $k\in\Z_{\geq 0}$. 

ii). We claim that $x[k]\in \mathcal{D}_2$ for some $k\in\Z_{\geq 0}$ implies that $x[k+n(N+1)]\in \mathcal{D}_2$ for any $n\in\Z_{\geq 0}$. Indeed, $x[k]\in \mathcal{D}_2$ implies that $V(x[k])\leq c_2+c_4$. If $x[k+N+1]\notin \mathcal{D}_2$, then $V(x[k+N+1])>c_2+c_4$. From \eqref{quasidetec2} we have $V(x[k])\geq V(x[k+N+1])-c_2>c_4$, which implies 
that $\eta_2p(x[k])>c_2+c_5$ by the definition of $c_4$. Then, again from  \eqref{quasidetec2}, we have $V(x[k+N+1])\leq V(x[k])+c_2 -\eta_2p(x[k])
< V(x[k]) \leq c_2+c_4$,
which contradicts with the assertion that $V(x[k+N+1])>c_2+c_4$. 
Define $j_s=\min\{k\in\Z_{\geq 0}\mid k\equiv s\; (\mbox{mod}\; N+1), V(x[k])\in\mathcal{D}_2\}\leq \infty$ for $s=0,1,...,N$ where ``$\mbox{mod}$'' denotes the modulo. 
The claim above shows that $x[k]\in \mathcal{D}_2$ for any $k\geq j_s$ where $k\equiv s\; (\mbox{mod}\; N+1)$. For any $k+N+1<j_s$ where $k\equiv s\; (\mbox{mod}\; N+1)$, $k\in\Z_{\geq 0}$, \eqref{quasidetec2} implies that $V(x[k+N+1])-V(x[k])\leq c_2 -\eta_2p(x[k])$. Since $x[k]\notin \mathcal{D}_2$, $V(x[k])>c_2+c_4>c_4$, which implies that $\eta_2p(x[k])>c_2+c_5$. Therefore,  $V(x[k+N+1])-V(x[k])\leq -c_5$. Hence, $j_s\leq (V(x[0]-c_2-c_4))/c_5<\infty$. Choose $K=\max\{j_0,j_1,...,j_N\}$. Then,  $x[k]\in\mathcal{D}_2$ for any $k\geq K,k\in\Z_{\geq 0}$. This completes the proof.
\hfill$\Box$

\section{Proof of Lemma \ref{lem4}}\label{prooflem4}
The system $\S_d$ is $0$-step SD implies that $\vartheta|u[k]|^2+|y[k]|^2\geq p(x[k]),\forall k\in\Z_{\geq 0}$. 
As in Lemma \ref{lembound}, we choose $\l$ such that $0<\l<\rho$, and define $\eta_1=\frac{1}{4\l}-\nu> 0$, $\eta_2=\rho-\l>0$. 
Then we have
\begin{align*}
V(x[k\!+\!1])\!-\!V(x[k])
\leq& \eta_1|u[k]|^2\!-\!\eta_2|y[k]|^2\!+\!\delta\nonumber\\
\leq& \!-\!\eta_2p(x[k])\!+\!(\eta_1\!+\!\vartheta\eta_2)|u[k]|^2\!+\!\delta\\
\leq& \!-\!\eta_2c|x[k]|^\lambda\!+\!(\eta_1\!+\!\vartheta\eta_2)|u[k]|^2\!+\!\delta\\
\end{align*}
Define $\mathcal{K}_\infty$-functions $\alpha_1(s)=as^\lambda$, $\alpha_2(s)=bs^\lambda$, $\alpha_3(s)=\eta_2cs^\lambda$, a $\mathcal{K}$-function $\sigma(s)=(\eta_1+\vartheta\eta_2)s^2$ and $d_2=\delta$. Then $V$ is an ISpS-Lyapunov function of $\S_d$ by Def. 10.  The conclusion follows by Theorem 2.5 in \cite{lazar2008input}. \hfill$\Box$


\section{Proof of Theorem \ref{thm1}} \label{sec:prf2}

Let the inputs of $\S$ to be piecewise constant, that is, $u(t)= u[k]$ for any $t\in[k\tau,(k+1)\tau)$. Since $x[k]=x(k\tau)$, we have
\begin{align}
&|y(t)-y[k]|
=|h_1(x(t))-h_1(x(k\tau))|\nonumber\\
&\leq \int_{k\tau}^{(k+1)\tau}|\dot{h}_1(s)|\;\ud s
\leq \sqrt{\tau}\sqrt{\int_{k\tau}^{(k+1)\tau}|\dot{h}_1(s)|^2\;\ud s},\label{app2dis}
\end{align}
where the last inequality is from the Cauchy-Schwarz inequality. For any $k_0,k_1\in\Z_{\geq 0},k_0<k_1$ and for any admissible $u(t)$, 
\begin{align*}
&\left|\int_{k_0\tau }^{k_1\tau}u(t)^\top y(t)\;\ud t-\tau\sum_{k=k_0}^{k_1-1} u[k]^\top y[k]\right|\\
\leq & \sum_{k=k_0}^{k_1-1}\int_{k\tau}^{(k+1)\tau}|u[k]||y(t)-y[k]|\;\ud t \\
\leq & \tau\sqrt{\tau}\sum_{k=k_0}^{k_1-1}|u[k]|\sqrt{\int_{k\tau}^{(k+1)\tau}|\dot{h}_1(s)|^2\;\ud s}  \nonumber\\
\leq & \tau\sqrt{\tau}\sqrt{\sum_{k=k_0}^{k_1-1}|u[k]|^2}\sqrt{\gamma^2\int_{k_0\tau}^{k_1\tau}|u(s)|^2\;\ud s+\beta(x[0])}  \nonumber\\
\leq & \tau^2\gamma \sum_{k=k_0}^{k_1-1}|u[k]|^2 + \frac{\tau\beta}{\gamma},
\end{align*}
where the second inequality is from \eqref{app2dis}, the third inequality is from the Cauchy-Schwartz inequality, the fourth inequality is from \eqref{eqngain}, the fifth inequality is because $\beta(x[k_0])>0$.
Therefore, 
\begin{align}
\int_{k_0\tau}^{k_1\tau}\!u(t)^\top \!y(t)\ud t&\leq \tau\!\sum_{k=k_0}^{k_1-1}  u[k]^\top\! y[k]\!+\!\tau^2\gamma \!\sum_{k=k_0}^{k_1-1}\!|u[k]|^2 \!+\! \frac{\tau\beta}{\gamma}.\label{Di1}
\end{align}
It is clear that
\begin{equation}\label{Di2}
-\nu\int_{k_0\tau}^{k_1\tau}|u(t)|^2\;\ud t=-\tau\nu\sum_{k=k_0}^{k_1-1} |u[k]|^2.
\end{equation}
Furthermore, 
\begin{align*}
& \left|\int_{k_0\tau}^{k_1\tau}y(t)^\top y(t)\;\ud t-\tau\sum_{k=k_0}^{k_1-1}y[k]^\top y[k]\right|\\
\leq& \sum_{k=k_0}^{k_1-1}\int_{k\tau}^{(k+1)\tau}|y(t)-y[k]|^2\;\ud t\\
&\quad+2\sum_{k=k_0}^{k_1-1}\int_{k\tau}^{(k+1)\tau}|y[k]||y(t)-y[k]|\;\ud t\\
\leq&\sum_{k=k_0}^{k_1-1}\int_{k\tau}^{(k+1)\tau}|y(t)- y[k]|^2\;\ud t\\
&\quad+\sum_{k=k_0}^{k_1-1}\int_{k\tau}^{(k+1)\tau}\l_1|y(t)- y[k]|^2
+\frac{1}{\lambda_1}| y[k]|^2\;\ud t\\
\leq&\tau^2(1\!+\!\l_1)\left(\gamma^2\!\!\int_{k_0\tau}^{k_1\tau}|u(t)|^2\;\ud t\!+\!\beta\right)\!+\!\frac{\tau}{\l_1}\sum_{k=k_0}^{k_1-1}| y[k]|^2\\
=&\tau^3\gamma^2 (1\!+\!\l_1)\sum_{k=k_0}^{k_1-1}|u[k]|^2\!+\!\frac{\tau}{\l_1}\sum_{k=k_0}^{k_1-1}|y[k]|^2\!+\!\tau^2 (1\!+\!\l_1)\beta
\end{align*}
where $\lambda_1$ is an arbitrary positive number, 
the third inequality is from \eqref{eqngain}. 
Therefore,  
\begin{align}
&-\rho\int_{k_0\tau}^{k_1\tau}y(t)^\top y(t)\;\ud t
\leq\tau^3\gamma^2 (1+\l_1)|\rho|\sum_{k=k_0}^{k_1-1}|u[k]|^2+\nonumber\\
&\;\tau(\frac{|\rho|}{\l_1}-\rho)\sum_{k=k_0}^{k_1-1}|y[k]|^2+\tau^2|\rho|(1+\l_1)\beta(x[k_0]).\label{Di3}
\end{align}
Since $V(x[k_0\tau])=V(x(k_0\tau))$, $V(x[k_1])=V(x(k_1\tau))$, $V(x(k_1\tau))-V(x(k_0\tau))
\leq\int_{k_0\tau}^{k_1\tau}u(t)^\top y(t)-\nu |u(t)|^2 - \rho |y(t)|^2\,\ud t$, 
adding (\ref{Di1}), (\ref{Di2}), (\ref{Di3}) results in \eqref{thm1:ineq1}. \hfill $\Box$
\section{Proof of Theorem \ref{thm2}} \label{sec:pfdegrade}
For any $k \in\Z_{\geq 0}$, since $|y[k]-\tilde y[k]|\leq \sqrt{m}\mu_2$, we have
\begin{align*}
\left||y[k]|^2-|\tilde y[k]|^2\right|\leq & 2|\tilde y[k]||y[k]-\tilde y[k]|\!+\!|y[k]-\tilde y[k]|^2\\
\leq & \frac{1}{\l_2}|\tilde y[k]|^2\!+\!(1+\l_2)m\mu_2^2
\end{align*}
where $\l_2$ is an arbitrary positive number. Hence,
\begin{align}
-\!\rho_1|y[k]|^2&\leq (\frac{|\rho_1|}{\l_2}\!-\!\rho_1)|\tilde y[k]|^2 \!+\! |\rho_1|(1\!+\!\l_2)m\mu_2^2.\label{Di13}
\end{align}
Similarly, since $|u[k]-\tilde u[k]|\leq \sqrt{m}\mu_1$, we have
\begin{align*}
\left||u[k]|^2-|\tilde u[k]|^2\right|\leq & |u[k]-\tilde u[k]|^2\!+\!2|\tilde u[k]||u[k]-\tilde u[k]|\\
\leq & \frac{1}{\l_3}|\tilde u[k]|^2\!+\!(1+\l_3)m\mu_1^2
\end{align*}
where $\l_3$ is an arbitrary positive number. Hence, 
\begin{align}\label{hatineq2}
-\!\nu_1|u[k]|^2&\leq ( \frac{|\nu_1|}{\l_3}\!-\!\nu_1 )|\tilde u[k]|^2 \!+\! |\nu_1|(1\!+\!\l_3)m\mu_1^2.
\end{align}
Since
\begin{align*}
&|u[k]^\top y[k]-\tilde u[k]^\top \tilde y[k]|\\
\leq &|\tilde u[k]||y[k]-\tilde y[k]|+|\tilde y[k]||u[k]-\tilde u[k]|\\
\leq &\l_4 m\mu_2^2 +\frac{1}{4\l_4}|\tilde u[k]|^2+ \l_5m\mu_1^2+\frac{1}{4\l_5}|\tilde y[k]|^2,
\end{align*}
where $\l_4,\l_5$ are arbitrary positive numbers, we have
\begin{align}\label{hatineq1}
u[k]^\top y[k]\leq&\tilde u[k]^\top \tilde y[k]+\frac{1}{4\l_4}|\tilde u[k]|^2+\frac{1}{4\l_5}|\tilde y[k]|^2\nonumber\\
&\quad +m(\l_4 \mu_2^2+\l_5\mu_1^2).
\end{align}
Then the inequality \eqref{thm2:ineq1} follows immediately from \eqref{thm1:ineq1}, \eqref{Di13}, \eqref{hatineq2} and \eqref{hatineq1}. \hfill $\Box$

\section{Proof of Theorem \ref{thm3}}\label{proofthm3}
Since $|\tilde u_2[k]-u_2[k]|\leq \sqrt{m}\mu_1$, $|\tilde u_2[k]|^2\geq \frac{1}{2}|u_2[k]|^2\!-\!|\tilde u_2[k]\!-\!u_2[k]|^2\geq\frac{1}{2}|u_2[k]|^2\!-\!m\mu_1^2$;
similarly, $|\tilde y_2[k]|^2\geq\frac{1}{2}|y_2[k]|^2-m\mu_2^2$. Hence, 
\begin{align}
&\sum_{k=k_0}^{k_0+N_2}\vartheta_2|\tilde u_2[k]|^2+|\tilde y_2[k]|^2\nonumber\\
\geq & \frac{1}{2}p_2(x_2[k_0])-m(N_2+1)(\vartheta_2\mu_1^2+\mu_2^2).\label{ineqdet2}
\end{align}
Then, similar to the proof of Lemma \ref{lemstrongdet}, it is easy to show that  $\sum_{k=k_0}^{k_0+N} \vartheta |r[k]|^2+|y[k]|^2\geq p(x[k_0])-d_2$. 

By Theorem \ref{thmfeedback}, the feedback-connected system satisfies inequality  \eqref{con:ineq3} where
$V(x)$ is positive definite, radially unbounded. Hence,
\begin{align}
& V(x[k_0+N+1])-V(x[k_0])\nonumber\\
\leq& \delta_1(x[k_0])\!+\!\delta_2(x[k_0])\!+\!\sum_{k=k_0}^{k_0+N} \left(\eta_1|r[k]|^2\!-\!\eta_2|y[k]|^2\!+\!\tilde \delta_2\right)\nonumber\\
\leq& d_1+d_2 -\kappa(x[k_0]).\label{ineqN}
\end{align}
Define $j_s=\min\{k\in\Z_{\geq 0}\mid k\equiv s\; (\mbox{mod}\; N+1), V(x[k])\in\mathcal{D}_2\}\leq \infty$ for $s=0,1,...,N-1$. As in part 2) of the proof of Lemma \ref{lembound}, 
it can be shown that 
$x[k]\in \mathcal{D}_2$ for any $k\geq j_s$ where $k\equiv s\; (\mbox{mod}\; N+1)$, and   $V(x[k_0+N+1])-V(x[k_0])\leq -d_3<0$ for any $k_0+N<j_s$ where $k_0\equiv s\; (\mbox{mod}\; N+1)$, $k_0\in\Z_{\geq 0}$. Then the conclusion follows immediately.
\hfill$\Box$

\section{Proof of Corollary \ref{cor4}}\label{proofthm4}
When $r[k]\equiv{\bf 0}$, $u_1[k]=-\tilde y_2[k]$ and $\tilde u_2[k]=y_1[k]$. Recalling \eqref{det1} and \eqref{ineqdet2}, we have $\sum_{k=k_0}^{k_0+N_1}\vartheta_1|\tilde y_2[k]|^2+|y_1[k]|^2\geq p_1(x_1[k_0])$,
$\sum_{k=k_0}^{k_0+N_2}\vartheta_2|y_1[k]|^2+|\tilde y_2[k]|^2\geq \frac{1}{2}p_2(x_2[k_0])
-m(N_2+1)(\vartheta_2\mu_1^2+\mu_2^2).$
Since $\tilde \nu_2+\rho_1'>0$,  $\nu_1'+\tilde\rho_2>0$, $\sum_{k=k_0}^{k_0+N} \big[(\tilde \nu_2+\rho_1')|y_1[k]|^2+( \nu_1'+\tilde\rho_2) |\tilde y_2[k]|^2\big]\geq p(x[k_0]) - d_2$. 
By Corollary \ref{corfeed}, when $r[k]\equiv{\bf 0}$, the feedback-connected system satisfies the inequality shown in \eqref{con:ineq5}. Then, $V(x[k_0+N+1])-V(x[k_0])\leq d_1+d_2-\kappa(x[k_0])$. 
The following proof is then the same as that of Theorem \ref{thm3} after obtaining \eqref{ineqN}.
\hfill $\Box$

\section{Proof of Theorem \ref{corsym}}\label{proofcor6}
Consider the system $T_\tau(P)\times_{\mathcal{F}}^w\tilde T_{\tau}(\S)$ shown in Fig. \ref{figfe3} where the disturbance $w[k]$ is chosen as 
\begin{align}
w[k]=\tilde y_2^s[k]-\tilde y_2[k],\;\forall k\in\Z_{\geq 0}.\label{eqw}
\end{align} 
Because $|y_2[k]-y_2^s[k]|= |h_{21}(x_2[k])-h_{21}(x_2^s[k])|\leq L|x_2[k]-x_2^s[k]|_\infty\leq L\epsilon$  by the Lipschitz property of $h_{21}$,  and $|\tilde y_2^s[k]-y_2^s[k]|\leq \sqrt{m}\mu_2$ by the property of the quantizer $\mathcal{Q}_2$, it follows that $|w[k]|\leq |y_2[k]-\tilde y_2[k]|+|y_2[k]-y_2^s[k]|+|y_2^s[k]-\tilde y_2^s[k]|\leq L\epsilon + 2\sqrt{m}\mu_2,\;\forall k\in\Z_{\geq 0}$. 
Hence, the particular choice of $w[k]$ shown in \eqref{eqw} satisfies assumption \eqref{wbound}.
Moreover, if both $T_\tau(P)\times_{\mathcal{F}}\tilde T_{\tau\mu\eta}(\S)$ and $T_\tau(P)\times_{\mathcal{F}}^w\tilde T_{\tau}(\S)$ are driven by the same $r[k]$, then $w[k]$ shown in \eqref{eqw} ensures that $\check y_2[k]=\tilde y_2^s[k]$ and $u_2[k]=u_2^s[k]$, $\forall k\in\Z_{\geq 0}$. By the definition of approximate bisimulation (cf. Def. \ref{dfnemubi}), $x_2[k]$ and $x_2^s[k]$ are related by $|x_2[k]- x_2^s[k]|_\infty\leq \epsilon,\;\forall k\in\Z_{\geq 0}$. 
The state boundedness of $T_\tau(P)\times_{\mathcal{F}}\tilde T_{\tau\mu\eta}(\S)$ can be derived immediately by Lemma \ref{thm5}. \hfill $\Box$

\end{spacing}
\end{document}